\documentclass[aop,citesort,MSNbibl,dvips]{arximspdf}
\usepackage{accents}
\usepackage{mathrsfs}
\usepackage{graphicx}

\doi{10.1214/10-AOP608}
\volume{39}
\issue{6}
\pubyear{2011}
\firstpage{2224}
\lastpage{2270}

\makeatletter
\newcommand{\dotcupdisplay}{\mbox{\,$\dot \cup$\,}}
\newcommand{\fraca}[2]{{#1/#2}}
\newcommand{\fracb}[2]{{(#1)/#2}}
\newcommand{\fracc}[2]{{#1/(#2)}}

\newcommand{\sun}{\mathbb{S}_{1}}

\newcommand{\R}{\mathbb R}
\newcommand{\N}{\mathbb N}
\newcommand{\C}{\mathbb C}

\newcommand{\op}[1]{\operatorname{#1}}
\newcommand{\tends}[2]{\overset{#1}{\underset{#2}{\longrightarrow}}}
\newcommand{\underset}[2]{\mathop{#2}\limits_{#1}}
\newcommand{\overset}{\stackrel}
\newcommand{\iint}{\int\!\!\int}
\newcommand{\notag}{\nonumber}
\newcommand{\cal}{\mathcal}

\newtheorem{theo}{Theorem}[section]
\newtheorem{proposition}[theo]{Proposition}
\newtheorem{lemma}[theo]{Lemma}
\newtheorem{cor}[theo]{Corollary}
\newproclaim{defi}[theo]{Definition}
\newproclaim{conjecture}{Conjecture}
\newproclaim{rek}[theo]{Remark}
\makeatother

\begin{document}
\begin{frontmatter}

\title{Random recursive triangulations of the disk via fragmentation theory}
\runtitle{Random recursive triangulations}

\begin{aug}
\author[A]{\fnms{Nicolas} \snm{Curien}\ead[label=e1]{nicolas.curien@ens.fr}}
\and
\author[B]{\fnms{Jean-Fran\c{c}ois} \snm{Le Gall}\corref{}\ead[label=e2]{jean-francois.legall@math.u-psud.fr}}
\runauthor{N. Curien and J.-F. Le Gall}
\affiliation{Ecole normale sup\'{e}rieure and
Universit\'{e} Paris-Sud}
\address[A]{DMA---Ecole normale sup\'{e}rieure\\
45 rue d'Ulm\\
75230 Paris Cedex 05\\
France\\
\printead{e1}} 
\address[B]{Math\'{e}matiques, b\^{a}t. 425\\
Universit\'{e} Paris-Sud\\
91405 Orsay Cedex \\
France\\
\printead{e2}}
\end{aug}

\received{\smonth{6} \syear{2010}}
\revised{\smonth{9} \syear{2010}}

%
\begin{abstract}
We introduce and study an infinite random triangulation of the unit
disk that
arises as the limit of several recursive models. This triangulation is
generated by throwing chords
uniformly at random in the unit disk and keeping only those chords that
do not intersect the
previous ones. After throwing infinitely many chords and taking the
closure of the resulting set, one gets
a random compact subset of the unit disk whose complement is a
countable union of triangles.
We show that this limiting random set has Hausdorff dimension $\beta^*
+1$, where $\beta^*=(\sqrt{17}-3)/2$, and
that it can be described as the
geodesic lamination coded by a random continuous function which is H\"
{o}lder continuous with
exponent $\beta^*-\varepsilon$, for every $\varepsilon>0$. We also
discuss recursive constructions
of triangulations
of the $n$-gon that give rise to the same continuous limit when $n$
tends to infinity.
\end{abstract}

%

\end{frontmatter}

\section{Introduction}\label{sec1}

\begin{figure}

\includegraphics{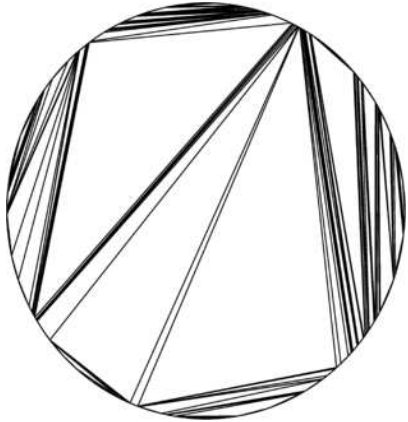}

\caption{The random set $L_\infty$.}
\label{fig1}
\end{figure}

In this work, we use fragmentation theory to study an infinite random
triangulation of the unit disk that arises as the limit of several
recursive models. Let us describe a special case of these models in
order to
introduce our main object of interest. We consider a sequence~$U_1,V_1,U_2,\allowbreak V_2,\ldots$ of independent
random variables, which are uniformly distributed over the unit circle
$\sun$
of the complex plane. We then construct inductively
a~sequence $L_1,L_2,\ldots$ of random closed subsets of the (closed)
unit disk $\overline{\mathbb D}$.
To begin with, $L_1$ just consists of the chord with endpoints $U_1$,
and $V_1$, which we denote
by $[U_1V_1]$. Then at step $n+1$, we consider two cases. Either the
chord $[U_{n+1}V_{n+1}]$
intersects $L_n$, and we put $L_{n+1}=L_n$, or the chord
$[U_{n+1}V_{n+1}]$ does not
intersect $L_n$, and we put $L_{n+1}=L_n \cup[U_{n+1}V_{n+1}]$. Thus,
for every integer $n\geq1$,
$L_n$ is a disjoint union of random chords. We then let
\[
L_\infty= \overline{\bigcup_{n=1}^\infty L_n}
\]
be the closure of the (increasing) union of the sets $L_n$.
See Figure~\ref{fig1} below for a simulation of the set $L_\infty$.

%

The closed set $L_\infty$ is a geodesic lamination of the unit disk,
in the sense that it
is a closed union of noncrossing chords (here we say that two chords
do not cross if they
do not intersect except possibly at their endpoints). We refer to \cite
{Bon98} for the general
notion of a geodesic lamination of a surface in the setting of
hyperbolic geometry. We may also view
$L_\infty$ as an infinite triangulation of the unit disk, in the same
sense as in Aldous \cite{Ald94b}.
Precisely,~$L_\infty$ is a closed subset of $\overline{\mathbb D}$,
which has zero Lebesgue measure
and is such that any connected component of $\overline{\mathbb
D}\setminus L_\infty$ is a triangle whose vertices
belong to the circle $\sun$.
The latter properties are not immediate, but will follow from
forthcoming statements.

In order to state our first result, let us introduce some notation. We~\mbox{denote}
the number of chords in $L_n$ by $N(L_n)$. Then, for every $x,y\!\in\!\sun$, we let~$H_{\!n}(x, y)$ be the number
of chords in $L_n$ that intersect the chord $[xy]$. We also set
\[
\beta^*=\frac{\sqrt{17}-3}{2}.
\]

\begin{theo}
\label{asympto}
\textup{(i)} We have
\[
n^{-1/2}N(L_n)\overset{a.s.}{\underset{n\to\infty}{\longrightarrow
}}\sqrt{\pi}.
\]

\textup{(ii)}
There exists a random process $(\mathscr{M}_{\infty}(x),x \in\sun
)$, which is H\"{o}lder continuous
with exponent $\beta^*-\varepsilon$, for every $\varepsilon>0$,
such that, for every $x\in\sun$,
\[
n^{-\beta^*/2} H_{n}(1,x) \overset{(\mathbb{P})}{\underset{n\to
\infty}{\longrightarrow}} \mathscr{M}_{\infty}(x),
\]
%
where $\overset{(\mathbb{P})}{\longrightarrow}$ denotes convergence
in probability.
\end{theo}

Part (i) of the theorem is a rather simple consequence of the results
in~\cite{BD86,BD87}, but
part (ii) is more delicate and requires different tools. In
the present work, we prove (more general versions of) the convergences
in (i) and (ii)
by using fragmentation theory. To this end, we consider continuous-time
models where
noncrossing
chords are thrown at random in the unit disk according to the following
device: At time~$t$, the existing chords bound several subdomains
of the disk, and a new chord is created in one of these subdomains at a
rate which is a given power of the
Lebesgue measure of the portion of the
circle that is adjacent to this subdomain. It is not hard to see that
the random closed subset of $\overline{\mathbb D}$
obtained by taking the closure of all chords created in this process
has the same distribution
as $L_\infty$, and moreover the case when the power is the square is
very closely related
to the discrete-time model described above.

In this continuous-time model, the ranked sequence of the Lebesgue
measures of the portions of $\sun$
corresponding to the subdomains bounded by the existing chords at time
$t$ forms
a conservative fragmentation process, in the sense of \cite{Ber06}. A
general version
of the convergence (i) can then be obtained
as a consequence of asymptotics for fragmentation processes.
Similarly, if $U$ is a random point uniformly distributed on $\sun$
and if we look
only at subdomains that intersect the chord $[1U]$, we get another
(dissipative) fragmentation process,
and known asymptotics give the convergence in (ii), provided that $x$
is replaced by
the random point~$U$. An extra absolute continuity argument is then
needed to get the desired result for a deterministic
point $x$: See Theorem \ref{tech} and its proof. It is plausible that
the convergence in (ii) also
holds almost surely, but the known asymptotics for fragmentation
processes do not give this stronger form.

The most technical part of the proof of Theorem \ref{asympto} is the
derivation of the
H\"{o}lder continuity properties of the limiting process $(\mathscr
{M}_\infty(x),x\in\sun)$. To this end,
we need to obtain precise bounds for the moments of increments of this
process. In order to
derive these bounds, we rely on integral equations for the moments,
which follow from
the recursive construction.

Our second theorem shows that the random geodesic lamination $L_\infty$
is coded by the process $\mathscr{M}_\infty$, in the sense of the
following statement.
For every $x,y\in\sun\setminus\{1\}$, we let $\operatorname{Arc}(1,x)=  \operatorname{Arc} (x,1)$ be the closed subarc of $\sun$ with endpoints
$x$ and $y$ that does not contain the point $1$. For every $x\in\sun
\setminus\{1\}$, we let
$\op{Arc}(1,x)=  \operatorname{Arc} (x,1)$ be the closed subarc of $\sun$ going from $1$ to $x$
in counterclockwise order,
and we set $\op{Arc}(1,1)=\{1\}$ by convention.

\begin{theo}
\label{codingL}
The following properties hold almost surely. The random set $L_\infty$
is the union of the chords $[xy]$ for all $x,y\in\sun$ such that
%
%
\begin{equation}
\label{maxipro}
\mathscr{M}_\infty(x)=\mathscr{M}_\infty(y)=\min_{z\in\op
{Arc}(x,y)} \mathscr{M}_\infty(z) .
\end{equation}
Moreover, $L_\infty$ is maximal for the inclusion relation among
geodesic laminations.\vadjust{\goodbreak}
\end{theo}

It is relatively easy to see that property (\ref{maxipro}) holds for
any chord $[xy]$ that
arises in our construction of $L_\infty$. The difficult part of the
proof of the theorem
is to show the converse, namely that any chord $[xy]$ such that (\ref
{maxipro}) holds
will be contained in $L_\infty$. This fact is indeed closely related
to the maximality
property of $L_\infty$.

The coding of geodesic laminations by continuous functions is discussed
in~\cite{LGP08}, and
is closely related to the coding of $\mathbb{R}$-trees by continuous functions
(see, e.g.,~\cite{DLG05}).
A particular instance of this coding had been discussed earlier by
Aldous \cite{Ald94b},
who considered the case when the coding function is the normalized
Brownian excursion. In that case, the associated $\mathbb{R}$-tree is
Aldous's CRT. Moreover, the Hausdorff dimension of the corresponding
lamination is $3/2$. This may be compared to the following statement, where
$\op{dim}(A)$ stands for the Hausdorff dimension of a subset $A$ of
the plane.

\begin{theo}
\label{Hausdim}
We have almost surely
\[
\op{dim}(L_\infty) = \beta^*+1 =\frac{\sqrt{17}-1}{2}.
\]
\end{theo}

The lower bound $\op{dim}(L_\infty)\geq\beta^*+1$ is a relatively
easy consequence
of the fact that $L_\infty$ is coded by the function ${\mathscr
M}_\infty$ (Theorem \ref{codingL}) and of the
H\"{o}lder continuity properties of this function (Theorem \ref
{asympto}). In order to get the
corresponding upper bound, we use explicit coverings of the set
$L_\infty$ that follow from our
recursive construction. To evaluate the sum of the diameters of balls
in these coverings
raised to a suitable power, we again use certain asymptotics from
fragmentation theory.

The random set $L_\infty$ also
occurs as the limit in distribution of certain random recursive
triangulations of the
$n$-gon. For every $n\geq3$, we consider the $n$-gon whose vertices
are the
$n$th roots of unity
\[
x^n_k=\exp\biggl(2i\pi\frac{k}{n}\biggr) , \qquad k=1,2,\ldots,n.
\]
A chord of $\sun$ is called a diagonal of the $n$-gon if its vertices
belong to the set $\{x^n_k\dvtx1\leq k\leq n\}$ and if it is not an
edge of
the $n$-gon.
A triangulation of the $n$-gon is the union of $n-3$ noncrossing
diagonals of the $n$-gon
(then the connected components of the complement of this union
in the $n$-gon are indeed triangles). The set $\mathscr{T}_{n}$ of all
triangulations of
the $n$-gon is in one-to-one correspondence with the set of all planar
binary trees with $n-1$ leaves (see, e.g.,
Aldous \cite{Ald94b}).

For every fixed integer $n\geq4$, we construct a random element of
$\mathscr{T}_{n}$ as
follows. Denote by $\mathscr{D}_n$ the set of all diagonals of the $n$-gon.
Let $c_1$ be chosen uniformly at random in $\mathscr{D}_n$. Then,
conditionally given $c_1$, let $c_2$ be a chord chosen uniformly at
random in the set
of all chords in $\mathscr{D}_n$ that do not cross~$c_1$. We
continue\vadjust{\goodbreak}
by induction and construct a finite
sequence of chords $c_1,c_2,\ldots,\allowbreak c_{n-3}$: for every $1<k\leq n-3$,
$c_k$ is chosen uniformly at
random in the set of all chords in $\mathscr{D}_n$ that do not cross
$c_1,c_2,\ldots,c_{k-1}$.
Finally we let~$\Lambda_n$ be the union of the chords $c_1,c_2,\ldots
,c_{n-3}$.

Let us also introduce a slightly different model, which is closely
related to~\cite{DHW08}. Let $\sigma$ be a
uniformly distributed random
permutation of $\{1,2,\ldots,n\}$. With $\sigma$, we associate a collection
of diagonals of the $n$-gon, which is constructed recursively as
follows. For every
integer $0\leq k\leq n$ we define a set $M_k$ of disjoint diagonals of
the $n$-gon, and a
set $F_k$ of ``free'' vertices. We start with $M_0=F_0=\varnothing$. Then,
at step $k\in\{1,\ldots,n\}$, either there is a (necessarily unique)
free vertex $x\in F_{k-1}$ such that
$[xx^n_{\sigma(k)}]$ is a diagonal of the $n$-gon that does not
intersect the chords in $M_{k-1}$, and we set
$M_k=M_{k-1}\cup\{[xx^n_{\sigma(k)}]\}$ and $F_k=F_{k-1}\setminus\{
x\}$; or there is no such vertex
and we set $M_k=M_{k-1}$ and $F_k =F_{k-1}\cup\{x^n_{\sigma(k)}\}$.
We let $\widetilde\Lambda_n$ be the union of the chords in $M_n$ (note that~$\widetilde\Lambda_n$
is \textit{not} a triangulation of the $n$-gon).

\begin{theo}
\label{discreteapprox}
We have
\[
\Lambda_{n} \tends{(d)}{n\to\infty} L_{\infty}
\]
and
\[
\widetilde\Lambda_{n} \tends{(d)}{n\to\infty} L_{\infty}.
\]
In both cases, the convergence holds in distribution in the sense of
the Hausdorff distance between compact subsets of $\overline{\mathbb D}$.
\end{theo}

Theorem \ref{discreteapprox} should be compared with the results of
Aldous \cite{Ald94b}
(see also~\cite{Ald94a}). Aldous considers a triangulation
of the $n$-gon that is uniformly distributed over $\mathscr{T}_n$, and
then proves that this random triangulation converges
in distribution
as $n\to\infty$ toward the geodesic lamination coded by the
normalized Brownian
excursion (see Theorem \ref{AldousTh} below for a more precise
statement). Our random recursive
constructions give rise to a limiting geodesic lamination which is
``bigger'' than the one that
appears in Aldous's work, in the sense of Hausdorff dimension.

Triangulations of convex polygons are also interesting from the
geometric and combinatorial point of view (see, e.g., \cite{STT88}).
In \cite{DFHN99}, Devroye et al. study some features of triangulations
sampled uniformly from $\mathscr{T}_{n}$. Their proofs are based on
combinatorial and enumeration techniques. Recursive triangulations of
the type studied in the present work have been used in physics as
greedy algorithms for computing folding of RNA structure (see \cite{Mul03}).
In these models, the polymer is represented by a discrete cycle and
diagonals correspond to liaisons of RNA bases. See \cite{Mul03,DHW08,DDJS09}
for certain results related to our work, and in
particular to the asymptotics of Theorem \ref{asympto}.

As a final remark, this work deals with ``Euclidean'' geodesic
laminations consisting of unions of chords.\vadjust{\goodbreak}
As in \cite{LGP08}, we may consider instead the
hyperbolic geodesic laminations obtained by replacing each chord by the
hyperbolic line
with the same endpoints in the hyperbolic disk. It is immediate to
verify that our main results remain valid after this replacement.

The paper is organized as follows. Section~\ref{sec2} recalls basic
facts about
geodesic laminations, and introduces the
random processes $(S_\alpha(t))_{t\geq0}$ describing random recursive
laminations, which are of interest in this work. Section~\ref{sec3}
studies the
connections between
these random processes and fragmentation theory, and derives general
forms of the asymptotics of Theorem~\ref{asympto}.
Section~\ref{sec4} is devoted to the continuity properties of the process~$\mathscr{M}_\infty$. Theorem \ref{codingL}
characterizing $L_\infty$ as the lamination coded by $\mathscr
{M}_{\infty}$ is proved in Section~\ref{sec5}. The
Hausdorff dimension of $L_\infty$ is computed in Section~\ref{sec6}, and
Section~\ref{sec7} discusses the discrete models
of Theorem~\ref{discreteapprox}. Finally, Section~\ref{sec8} gives some
extensions and comments.

\section{Random geodesic laminations}\label{sec2}
\subsection{Laminations}
\label{laminations}
Let us briefly recall the notation which was already introduced in
Section~\ref{sec1}.
The open unit disk of the complex plane $\C$ is denoted by $\mathbb
{D}=\{z\in\mathbb{C} \dvtx |z|<1\}$, and $\sun$ is the unit circle.
As usual, the closed unit disk is denoted by $\overline{\mathbb{D}}$.
If $x,y$ are two distinct
points of $\sun$, the \textit{chord of feet} $x$ \textit{and} $y$ is
the closed line segment $[xy]\subset\overline{\mathbb{D}}$. We also
use the
notation $]xy[$ for the open line segment with endpoints $x$ and $y$.
By convention, $[xx]$ is equal to
the singleton $\{x\}$, and is viewed as a degenerate chord, with
$]xx[=\varnothing$.

We say that two chords $[xy]$
and $[x'y']$ do not cross if
$]xy[\,\cap\,]x'y'[\,=\varnothing$.

\begin{defi}
A geodesic lamination $L$ of $\overline{\mathbb{D}}$ is a closed
subset $L$
of $\overline{\mathbb{D}}$ which can be written as
the union of a collection of noncrossing chords.
The lamination $L$ is \textit{maximal} if it is maximal for the
inclusion relation among geodesic laminations of
$\overline{\mathbb{D}}$.
\end{defi}

For simplicity, we will often say lamination instead of geodesic
lamination of~$\overline{\mathbb{D}}$.
In the context of hyperbolic geometry \cite{Bon98}, geodesic
laminations of the disk are defined as closed subsets of the open
(hyperbolic) disk.
Here we prefer to view them as compact subsets of the closed disk,
mainly because we want to discuss convergence of laminations in the sense
of the Hausdorff distance. Notice that a maximal lamination necessarily
contains the unit circle $\sun$.

As the next lemma shows, the concept of a maximal lamination is a
continuous analogue of a discrete triangulation.

\begin{lemma}\label{triangle} Let $L$ be a geodesic lamination of
$\overline{\mathbb{D}}$. Then $L$ is maximal if and only if the
connected components of $\overline{\mathbb{D}}\setminus L$ are open
triangles whose vertices belong to $\sun$.
\end{lemma}

We leave the easy proof to\vadjust{\goodbreak} the reader.

\subsection{Figelas and associated trees}\label{sec2.2}
The simplest examples of laminations are finite unions of noncrossing
chords. Define a \textit{figela} $S$ (from \textit{fi}nite \textit
{ge}odesic \textit{la}mination) as a finite set of (unordered) pairs of
distinct points of $\sun$
\[
S=\{\{x_{1},y_{1}\},\ldots,\{x_{n},y_{n}\}\},
\]
such that the union of the $n$ chords $\{[x_{i}y_{i}]\}_{1\leq i \leq
n}$ forms a lamination, which~is then denoted by $L_{S}$. If $\{x,y\}\,{\in}\,S$,
we will say that $[xy]$ is a chord of the fige\-la~$S$. We denote
the set $\bigcup_{i=1}^n\{x_{i},y_{i}\}$ of all feet of the chords of
$S$ by $\op{Feet}(S)$.

Let $u,v \in\sun\setminus\op{Feet}(S)$. The \textit{height}
between $u$ and $v$ in $S$ is the number of chords of $S$ crossed by
the chord $[uv]$
\[
H_{S}(u,v)=\#\{ 1\leq i \leq n \dvtx [x_{i}y_{i}]\cap[uv]\ne
\varnothing\}.
\]
The next proposition follows from simple geometric considerations.

\begin{proposition}[(Triangle inequality)] \label{ineg} Let $S$ be a figela.
For every $x,y,z\in\sun\setminus\op{Feet}(S)$ we have
%
%
\begin{equation} H_{S}({x},{z}) \leq H_{S}({x},{y}) +H_{S}({y},{z}).
\end{equation}
\end{proposition}

\begin{figure}

\includegraphics{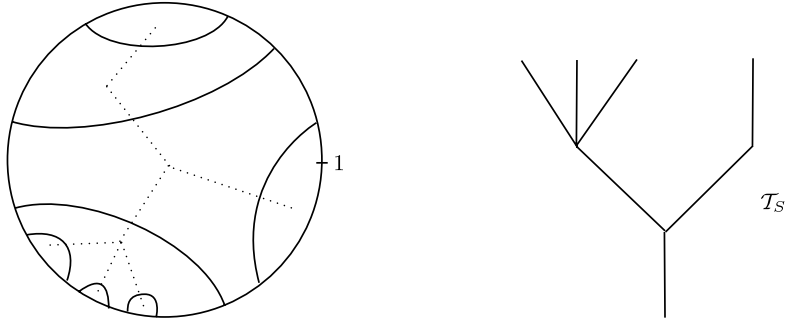}

\caption{A figela and its associated plane tree (in dotted lines on
the left side). We drew chords as
curved lines for better visibility. In this example, $S$ has $7$ chords
and $8$ fragments. Notice that each fragment of $S$ corresponds to a
vertex of the tree ${\mathcal T}_{S}$.}
\label{fig2}
\end{figure}

Let $S=\{\{x_{1},y_{1}\},\ldots,\{x_{n},y_{n}\}\}$ be a figela.
We define an equivalence relation on $\sun\setminus\op{Feet}(S)$ by
setting, for every $u,v\in\mathbb{S}_{1}\setminus\mbox{Feet}(S)$,
\[
u\simeq v \quad\mbox{if and only if} \quad H_{S}(u,v)=0.
\]
In other words, two points of $\sun\setminus\op{Feet}(S)$ are
equivalent if and only if they belong to the same connected component
of $\overline{\mathbb{D}} \setminus\bigcup_{i=1}^n[x_iy_i]$.
Then $H_S$ induces a distance on the quotient set ${\mathcal
T}_S:=(\sun\setminus\op{Feet}(S))/\!\! \simeq$.
The finite metric space~${\mathcal T}_{S}$ can be viewed as a graph by
declaring that there is an edge between~$a$ and $b$ if and only if
$H_{S}(a,b)=1$. This graph is indeed a tree, and~$H_{S}(\cdot,\cdot)$ coincides
with the usual graph distance. The tree
${\mathcal T}_{S}$ can be rooted at the equivalence class of $1$ (we
assume that $1$ is not a foot of $S$, which will always be the case
in our examples). As a result of this discussion, we can associate a
plane (rooted ordered) tree ${\mathcal T}_{S}$ to $S$. See Figure~\ref
{fig2} for
an example from which the definition of the tree ${\mathcal T}_S$
should be clear.

The $n+1$ connected components of $\overline{\mathbb{D}} \setminus
\bigcup_{i=1}^n[x_iy_i]$ are called the \textit{fragments} of the figela
$S$. With each fragment $R$, we associate
its \textit{mass}
\[
\mathrm{m}(R)= \lambda(R \cap\sun),
\]
where $\lambda$ denotes the uniform probability measure on
$\sun$.

%

\subsection{Coding by continuous functions} \label{coding}
Let $g\dvtx[0,1] \to\R_{+}$ be a continuous function such that
$g(0)=g(1)=0$. We define a pseudo-distance on $[0,1]$ by
\[
\mathrm{d}_{g}(s,t) = g(s)+g(t)-2\min_{r \in[s\wedge t,s\vee t]}
g(r)
\]
for every $s,t\in[0,1]$.
The associated equivalence relation on $[0,1]$ is defined\vspace*{-1pt} by setting
$s\overset{g}{\sim}t$ if and only if $\mathrm{d}_{g}(s,t)=0$, or
equivalently\vspace*{-2pt} $g(s)=g(t)=\min_{r \in[s\wedge t,s\vee t]} g(r)$.\vadjust{\goodbreak}

\begin{proposition}[(\cite{DLG05})] The quotient set $T_{g}:=[0,1]/\!\!
\overset{g}{\sim}$ endowed with the distance $\mathrm{d}_{g}$ is an $\R
$-tree called the tree coded by the function $g$.
\end{proposition}

We refer to \cite{Eva08} for an extensive discussion of $\R$-trees in
probability theory.

In order to introduce the lamination coded by $g$, we need some
additional notation. For $s\in[0,1]$,
we let $\op{cl}_g(s)$ be the equivalence class of $s$ with respect\vspace*{-2pt} to
the equivalence relation
$\overset{g}{\sim} $. Then, for $s,t\in[0,1]$, we set $s\overset
{g}{\approx} t$ if at least one of the following two conditions holds:
\begin{itemize}
\item$s\overset{g}{\sim} t$ and $g(r)> g(s)$ for every $r \in\,]s\wedge t,s\vee t[$.
\item$s\overset{g}{\sim} t$ and $s\wedge t =\min\op{cl}_g(s)$,
$s\vee t=\max\op{cl}_g(s)$.
\end{itemize}

In particular, $s\overset{g}{\approx}s $, and $s\overset{g}{\approx
} t$ holds if and only if
$t\overset{g}{\approx} s$. Note, however, that $\overset{g}{\approx
}$ is in general not an
equivalence relation. It is an elementary exercise to check that the
graph $\{(s,t)\dvtx s\overset{g}{\approx} t\}$
is a closed subset of $[0,1]^2$.

\begin{proposition} \label{codinglamination} The set
%
%
\begin{equation}
\label{lamicoding}
L_{g}:= \bigcup_{s\overset{g}{\approx}t} [e^{2i\pi s}e^{2i\pi
t} ]
\end{equation}
is a geodesic lamination of $\overline{\mathbb D}$ called the
lamination coded by the function $g$.
Furthermore, $L_g$ is maximal if and only if, for every open
subinterval $]s,t[$ of $[0,1]$, the
infimum of $g$ over $]s,t[$ is attained at at most one point of $]s,t[$.
\end{proposition}

We leave the proof to the reader. See \cite{LGP08}, Proposition~2.1,
for a closely related statement. This proposition is stated under the
assumption that the local minima
of $g$ are distinct, which is slightly stronger than the condition in
the second assertion of
Proposition \ref{codinglamination}. Note that the latter condition is
equivalent to saying that
the relations $\overset{g}{\sim}$ and $\overset{g}{\approx}$
coincide, or that $\overset{g}{\approx}$
is an equivalence relation.

We end this section by reformulating in this formalism a theorem of
Aldous which was already mentioned
in Section~\ref{sec1}. Recall our notation
$\mathscr{T}_{n}$ for the set of all triangulations of the $n$-gon. An
element of $\mathscr{T}_{n}$
is just a geodesic lamination consisting of $n-3$ chords whose feet
belong to the set of
$n$th roots of unity.

\begin{theo}[(\cite{Ald94b,Ald94a})]
\label{AldousTh}
Let $(\mathbf{e}_{t})_{0\leq t \leq1}$ be a normalized Brownian
excursion, and let $\Delta_{n}$ be uniformly distributed over
$\mathscr{T}_{n}$. Then we have
\[
\Delta_{n} \tends{(d)}{n\to\infty} L_{\mathbf{e}},
\]
in the sense of the Hausdorff distance between compact subsets of
$\mathbb{\overline{D}}$.
Moreover the Hausdorff dimension of $L_{\mathbf{e}}$ is almost surely
equal to $3/2$.
\end{theo}

A detailed argument for the calculation of the Hausdorff dimension of
$L_{\mathbf{e}}$
is given in \cite{LGP08} (the proof in \cite{Ald94b} is only sketched).

\subsection{Random recursive laminations}
\label{figelapro}
Let $\alpha\geq0$ be a positive real number. We define a Markov jump process
$(S_{\alpha}(t),t\geq0)$ taking values in the space of all figelas,
and increasing in the sense of the inclusion order.

Let us describe the construction of this process.
We introduce a sequence of random times $0=\tau_{0}<\tau_{1}< \tau
_{2}< \cdots$ such that $S_{\alpha}(t)$ is
constant over each interval $[\tau_{n},\tau_{n+1}[$, and $S_{\alpha
}(\tau_{n})$ has exactly $n$ chords [in particular $S_{\alpha}(0)$ is
the empty figela]. We define the pairs $(\tau_{n},S_{\alpha}(\tau_{n}))$
for every $n\geq1$
recursively as follows.
In order to describe the joint distribution of $(\tau_{n+1},S_{\alpha
}(\tau_{n+1}))$ given the $\sigma$-field $\mathcal{F}_{n}= \sigma
(\tau_{0}, \ldots , \tau_{n}, S_{\alpha}(\tau_{1}), \ldots ,
S_{\alpha
}(\tau_{n}))$, we write $R_{1}^n, \ldots , R_{n+1}^n$ for the $n+1$
fragments of the figela $S_\alpha(\tau_{n})$, and we let ${e}_{1},
\ldots , {e}_{n+1}$ be $n+1$ independent exponential variables with
parameter $1$ that are also independent of $\mathcal{F}_{n}$. Then,
for $1\leq j\leq n+1$, we set $\mathcal{E}_{j} = \mathrm
{m}(R_{j}^n)^{-\alpha}{e}_{j}$ and we let~$j_{0}$ be the a.s. unique
index such that $\mathcal{E}_{j_{0}}=\min\{\mathcal{E}_{j}\dvtx
1\leq
j\leq n+1\}$. Conditionally given $\mathcal{F}_{n}$ and $({e}_{1},
\ldots
, {e}_{n})$, we sample two independent random variables $X_{n+1},$ and
$Y_{n+1}$ uniformly distributed over $R_{j_{0}}\cap\sun$. Then
conditionally on $\mathcal{F}_{n}$, the pair $(\tau_{n+1},S_{\alpha
}(\tau_{n+1}))$ has the same distribution as $(\tau_{n}+\mathcal
{E}_{j_{0}}, S_{\alpha}(\tau_{n})\cup\{\{X_{n+1},Y_{n+1}\}\})$.

Note that $\tau_{n} \to\infty$ a.s. when $n \to\infty$. Indeed,
it is enough to see this when $\alpha=0$, and then
$\tau_{n+1}-\tau_n$ is exponential with parameter $n+1$. Therefore
the processes $S_\alpha(t)$ are well-defined
for every $t\geq0$.

If $R$ is a fragment of $S_{\alpha}(t)$, then independently of the
past up to time~$t$, a new chord is added in $R$ at rate $\mathrm
{m}(R)^{\alpha}$.
The preceding construction can thus be interpreted informally: the
first chord is thrown in $\overline{\mathbb{D}}$
uniformly at random (the two endpoints of the chord are chosen independently
and uniformly over $\sun$) at an exponential time with parameter $1$,
and divides it into two fragments $R_{0}$ and $R_{1}$. These two
fragments can be identified with two disks if we contract the first
chord (the boundaries of these disks are then identified,
respectively, with $R_{i}\cap\sun$ for $i\in\{0,1\}$). Then the
process goes on independently inside each of these disks provided that
we rescale time by the mass of the corresponding fragment to the power~$\alpha$.

The process $(S_{\alpha}(t),t\geq0)$ will be called the figela
process with autosimilarity parameter $\alpha$.

\begin{rek}
\label{Poissonconst}
Let $\{(t_{i},(x_{i},y_{i}))\}_{i \in\N}$ be the atoms of a
Poisson point measure on $\R_{+}\times\sun\times\sun$ with
intensity $dt \otimes\lambda\otimes\lambda$, where
we recall that $\lambda$ is the uniform probability measure on $\sun
$. We suppose that the atoms of the Poisson measure are ordered so that
$0<t_{1}<t_{2}<\cdots,$ and we also set
$t_0=0$. We construct a figela-valued jump process
$(\mathscr{S}(t),t\geq0)$ using the following device. We start from
$\mathscr{S}(0)=\varnothing$, and the process may jump only at times~$t_1,t_2,\ldots.$ For every
$i\geq1$, we take $\mathscr{S}(t_i)=\mathscr{S}(t_{i-1})\cup\{\{
x_i,y_i\}\}$
if the chord of feet $x_i$ and $y_i$ does not cross any chord
of $\mathscr{S}(t_{i-1})$, and otherwise we take $\mathscr
{S}(t_{i})=\mathscr{S}(t_{i-1})$. It follows from properties of
Poisson measures that this process has the same law as our process
$(S_{2}(t),t\geq0)$.
Moreover, the discrete-time process $(L_{\mathscr{S}(t_n)},n\geq0)$
has the same distribution
as the process $(L_n,n\geq0)$ discussed in Section~\ref{sec1}. Thanks to this
observation, and
to the fact that $n^{-1}t_n$ tends to $1$ a.s., forthcoming results
about asymptotics of the processes
$(S_{\alpha}(t),t\geq0)$ will carry over to the process $(L_n,n\geq
0)$.
\end{rek}

\begin{rek}[(Rotational invariance)] \label{rotation} Let
$(S_{\alpha}(t))_{t\geq0}$ be a figela process with parameter $\alpha
$. For every $z \in\sun$,
set
\[
{S}^z_{\alpha}(t) = \bigl\{\{zx,zy\} \dvtx \{x,y\} \in S_{\alpha
}(t) \bigr\}.
\]
Then $({S}^z_{\alpha}(t),t\geq0)$ has the same distribution as
$({S}_{\alpha}(t),t\geq0)$.
\end{rek}

It will be important to construct simultaneously the processes
$(S_\alpha(t))_{t\geq0}$ for all
values of $\alpha\geq0$, in the following way. We set
%
%
\begin{equation}
\label{binarytree}
{\mathbb T}=\bigcup_{n\geq0} \{0,1\}^n,
\end{equation}
where $\{0,1\}^0=\{\varnothing\}$. We
consider a collection $(\epsilon_u)_{u\in{\mathbb T}}$ of independent
exponential variables
with parameter $1$. The first chord then appears at time $\epsilon
_\varnothing$. If $R_0$ and $R_1$
are the two fragments created at this moment, a new chord will appear
in $R_0$, respectively, in
$R_1$, at time $\epsilon_\varnothing+ \mathrm{m}(R_0)^{-\alpha}\epsilon
_0$, respectively, at time
$\epsilon_\varnothing+ \mathrm{m}(R_1)^{-\alpha}\epsilon_1$. We
continue the construction by
induction. If we use the same random choices of the new chords
independently of
$\alpha$ (so that the same fragments will also appear), we get a
coupling of the processes
$(S_\alpha(t))_{t\geq0}$ for all $\alpha\geq0$.

This coupling is such that a.s. for every $t \geq0$ and for every
$\alpha' \geq\alpha\geq0$, there exists a finite random time
$T_{t,\alpha,\alpha'} \geq t$ such that
\[
S_{\alpha'}(t) \subset S_{\alpha}(t) \subset S_{\alpha'}(T_{t,\alpha
,\alpha'}).
\]

In the remaining part of this work, we will always assume that the
processes $(S_{\alpha}(t))_{t\geq0}$ are coupled in this way. Hence,
the increasing limit $S(\infty)=\lim\uparrow S_\alpha(t)$
as $t\uparrow\infty$ does not depend on $\alpha$, and the same holds
for the random closed subset of $\overline{\mathbb D}$
defined by
\[
L_{\infty} = \overline{\bigcup_{\{x,y\} \in S(\infty) }[xy]}.
\]
By the discussion in Remark \ref{Poissonconst}, this is consistent
with the definition of
$L_\infty$ in Section~\ref{sec1}. We note that
$L_\infty$ is a (random)
geodesic lamination. To see this,
write $S^*(\infty)$ for the closure in $\sun\times\sun$ of the set
of all (ordered) pairs $(x,y)$ such that $\{x,y\}$ belongs to $S(\infty
)$. Then a simple argument shows that
\[
L_{\infty} = \bigcup_{(x,y) \in S^*(\infty)} [xy],
\]
and moreover if $(x,y)$ and $(x',y')$ belong to $S^*(\infty)$ the
chords $[xy]$ and $[x'y']$ either coincide or
do not cross.

\section{Random fragmentations}\label{sec3}
\subsection{Fragmentation theory}
In this subsection, we briefly recall the results from fragmentation
theory that we will use, in the particular
case of binary fragmentation which is relevant to our applications. For
a more detailed
presentation, we refer to Bertoin's book \cite{Ber06}.

We
consider a probability measure $\nu$ on $[0,1]^2$. We assume that $\nu
$ is supported on the set $\{(s_{1},s_{2})\dvtx 1>s_{1}\geq s_{2}\geq0,
s_{1}+s_{2}\leq1\}$, and satisfies the following additional properties:
\renewcommand{\theequation}{H}
%
%
\begin{eqnarray}
\mbox{(i)} \quad\nu(\{s_{2}>0\})>0. \nonumber
\\[-8pt]
\\[-8pt]
\mbox{(ii)} \quad\nu(\{s_{1}=0\})=0.
\nonumber
\end{eqnarray}
\renewcommand{\theequation}{\arabic{equation}}
\setcounter{equation}{4}

Such a measure is a special case of a \textit{dislocation measure}.
Furthermore, if $\nu(\{s_{1}+s_{2}=1\})=1$, then $\nu$ is said to be
\textit{conservative}. It is called \textit{nonconservative} or
\textit{dissipative} otherwise.

Let $\mathcal{S}^{\downarrow}$ be the set of all real sequences
$(s_{1},s_{2}, \ldots)$ such that $1 \geq s_{1}\geq s_{2} \geq\cdots
\geq
0$ and $ \sum_{i=1}^\infty s_{i} \leq1$. A fragmentation process with
autosimilary parameter $\alpha\geq0$, and dislocation measure $\nu$
is a Markov process \mbox{$(X^{(\alpha)}(t),t\geq0)$} with values in
$\mathcal{S}^{\downarrow}$ whose evolution can be described
informally as follows
(see~\cite{Ber06} for a more rigorous presentation). Let $X^{(\alpha
)}(t) = (s_{1}(t),s_{2}(t), \ldots)$ be the state of the process at time
$t\geq0$. For each $i \geq1$,
$s_i(t)$ represents the mass of the $i$th particle at time $t$
(particles are ranked according to decreasing masses). Conditionally on
the past up to time $t$, the $i$th particle lives after time $t$ during
an exponential time of parameter $(s_{i}(t))^\alpha$, then dies and
gives birth to two particles of respective masses $R_{1}s_{i}(t)$ and
$R_{2}s_{i}(t)$, where the pair $(R_{1},R_{2})$ is sampled from $\nu$
independently of the past.

\begin{rek} We will not be interested in the case $\alpha<0$,
which is not relevant for our applications.
\end{rek}

We can construct
simultaneously the processes $(X^{(\alpha)}(t))_{t\geq0}$ starting
from $X^{(\alpha)}(0)=(1,0,\ldots)$, for all values of $\alpha\geq
0$ in the following way.
Consider first the process $X^{(0)}$ corresponding to $\alpha=0$.
We represent the genealogy of this process
by the infinite binary tree $\mathbb T$ defined in (\ref{binarytree}).
Each $u\in{\mathbb T}$ thus corresponds to a ``particle'' in the
fragmentation process. We denote the mass of $u$ by
$\xi_u$ and the lifetime of $u$ by $\zeta^{(0)}_u$. Since we are
considering the case
$\alpha=0$, the random variables $(\zeta^{(0)}_u)_{u\in\mathbb{T}}$
are independent and exponentially
distributed with parameter $1$. If we now want to construct
$(X^{(\alpha)}(t))_{t\geq0}$ for a given value
of $\alpha$, we keep the same values $\xi_u$ for the masses of
particles, but we replace
the lifetimes by $\zeta^{(\alpha)}_u=(\xi_u)^{-\alpha}\zeta
^{(0)}_u$, for every $u\in\mathbb{T}$. See
\cite{Ber06}, Corollary~1.2, for more details.

In the remaining part of this subsection, we assume that the processes
$(X^{(\alpha)}(t))_{t\geq0}$
starting from $X^{(\alpha)}(0)=(1,0,\ldots)$ are defined
for every $\alpha\geq0$ and coupled as explained above.

We set for every real $p \geq0$,
\[
\kappa_{\nu}(p) = \int_{[0,1]^2} \bigl(1-(s_{1}^p+s_{2}^p) \bigr)
\nu(ds_{1},ds_{2}),
\]
where by convention $0^0=0$. Then $\kappa_{\nu}$ is a continuous
increasing function. Under Assumption $(H)$, $\kappa_{\nu}(0)<0$ and
$\kappa_{\nu}(+\infty)=1$, and therefore there exists a unique
$p^*>0$, called the Malthusian exponent of $\nu$, such that
\[
\kappa_{\nu}(p^*)=0.
\]
The Malthusian exponent allows us to introduce the so-called Malthusian
martingale, which
is discussed in part (i) of the next theorem.

\begin{theo} \label{outil} Write $X^{(\alpha)}(t) = (s_{1}^{(\alpha
)}(t),s_{2}^{(\alpha)}(t), \ldots)$
for every $t\geq0$ and \mbox{$\alpha\geq0$}. Then:
\begin{longlist}[(iii)]
\item[(i)] For every $\alpha\geq0$, the process
\[
\mathscr{M}^{(\alpha)}(t):= \sum_{i=1}^\infty\bigl(s_{i}^{(\alpha
)}(t) \bigr)^{p^*},
\qquad t\geq0,
\]
is a uniformly integrable martingale and converges almost surely to a
limiting random variable $\mathscr{M}_{\infty}$, which does not
depend on $\alpha$.
Moreover $\mathscr{M}_{\infty}>0$ a.s., and $\mathscr{M}_{\infty}$
satisfies the following
identity in distribution:
%
%
\begin{equation}
\label{identdistri}
\mathscr{M}_{\infty}\overset{(d)}{=}
\Sigma_{1}^{p^*}\mathscr{M}_{\infty}' + \Sigma_{2}^{p^*}\mathscr
{M}_{\infty}'',
\end{equation}
where $(\Sigma_{1},\Sigma_{2})$ is distributed according to $\nu$,
and $\mathscr{M}_{\infty}'$ and $\mathscr{M}_{\infty}''$ are
independent copies of $\mathcal{\mathscr{M}_{\infty}}$, which are
also independent of the pair $(\Sigma_{1},\Sigma_{2})$. This identity
in distribution characterizes
the distribution of $\mathscr{M}_{\infty}$ among all probability
measures on ${\mathbb R}_+$ with mean $1$. Furthermore, we have
$\mathbb{E} [ \mathscr{M}_\infty^q ]<\infty$
for every real $q\geq1$.
\item[(ii)] For every real $p \geq0$, the process
\[
e^{t\kappa_{\nu}(p)}\sum_{i=1}^\infty\bigl(s^{(0)}_{i}(t)\bigr)^p, \qquad
t\geq0,
\]
is a martingale
and converges a.s. to a positive limiting random variable.
\item[(iii)] Let $\alpha>0$. Assume that $\int s_2^{-a}\nu
(ds_1,ds_2)<\infty$
for some $a>0$. Then for every $p\geq0$,
\[
t^{\fracb{p-p^*}{\alpha}}\sum_{i=1}^\infty\bigl(s^{(\alpha)}_{i}(t)\bigr)^p
\overset{\mathbb{L}^2}{\underset{t\to\infty}{\longrightarrow}}
K_{\nu}(\alpha,p)\mathscr{M}_{\infty},
\]
where $K_{\nu}(\alpha,p)$ is a positive constant depending on
$\alpha$, $p$ and $\nu$, and the limiting variable $\mathscr
{M}_{\infty}$ is the same
as in \textup{(i)}.\vspace*{-1pt}
\end{longlist}
\end{theo}

\begin{pf} The fact that $\mathscr{M}^{(\alpha)}(t)$
is a uniformly integrable martingale follows from \cite{Ber06},
Proposition~1.5. This statement also
shows that the almost sure limit $\mathscr{M}_{\infty}$
of this martingale coincides with the limit of the so-called intrinsinc
martingale, and therefore
does not depend on $\alpha$. By uniform integrability,
we have $\mathbb{E} [ \mathscr{M}_{\infty} ]={\mathbb E}[\mathscr
{M}^{(\alpha
)}(0)]=1$. The property $\mathscr{M}_{\infty}>0$ a.s.
follows from \cite{Ber06}, Theorem~1.1. The identity in distribution
(\ref{identdistri})
is a special case of (1.20) in \cite{Ber06}. The fact that the
distribution of
$\mathscr{M}_{\infty}$ is characterized by this identity (and the
property $\mathbb{E} [ \mathscr{M}_{\infty} ]=1$)
follows from Theorem~1.1 in \cite{Liu97}. The property $\mathbb{E} [
\mathscr {M}_\infty^q ]<\infty$ for every
$q\geq1$ is a consequence of Theorem~5.1 in the same article.

Then, assertion (ii) follows from Corollary~1.3 and Theorem~1.4 in
\cite{Ber06}.
Finally, assertion (iii) can be found in \cite{BG04}, Corollary~7,
under more general assumptions.\vspace*{-1pt}
\end{pf}

\begin{rek}
\label{conscase}
In the conservative case, we immediately
see that $p^*=1$ and $\mathscr{M}_\infty=1$.\vspace*{-1pt}
\end{rek}

\subsection{The number of chords in the figela process}

Let $\nu_C$ be the probability measure on $[0,1]^2$ defined by
\[
\int\nu_{C}(ds_{1},ds_{2}) F(s_{1},s_{2}) = 2\int_{1/2}^1 du
F(u,1-u)\vadjust{\goodbreak}
\]
for every nonnegative Borel function $F$. Clearly $\nu_C$ satisfies
the assumptions of the
previous subsection.\vspace*{-1pt}

\begin{proposition}\label{nbre} Fix $\alpha\geq0$. We denote by
$R_{1}^{\alpha}(t),R_{2}^{\alpha}(t),
\ldots$ the fragments of the figela $S_{\alpha}(t)$, ranked according
to decreasing masses. Then the process
\[
X_{{\alpha}}(t) = (\mathrm{m}(R_{1}^\alpha(t)),\mathrm{m}(R_{2}^\alpha(t)),
\ldots ),
\]
is a fragmentation process with parameters $(\alpha,\nu_{C})$.\vspace*{-1pt}
\end{proposition}

\begin{pf} From the construction of the figela processes, we see that, when
a chord appears in a fragment $R$ of the figela, it divides this
fragment into two new fragments of respective masses $U\mathrm{m}(R)$ and
$(1-U)\mathrm{m}(R)$ where $U$ is uniformly distributed over $[0,1]$. The ranked
pair of these masses is thus distributed as $(s_1\mathrm{m}(R),s_2\mathrm
{m}(R))$ under $\nu_C(ds_1,ds_2)$.
Furthermore a fragment~$R$ splits at rate $\mathrm{m}(R)^\alpha$. The
desired conclusion easily follows.
We leave details to the reader.\vspace*{-1pt}
\end{pf}

\begin{rek}
\label{couplingremark}
The coupling of $(S_{\alpha}(t),t\geq0)$ for all $\alpha\geq0$
yields a coupling of the associated fragmentation processes $(X_{\alpha
}(t),t\geq0)$. This is indeed the same coupling that was already
discussed in the previous subsection.\vspace*{-1pt}
\end{rek}

By combining Proposition \ref{nbre} with Theorem \ref{outil}, we
already get detailed information about the asymptotic number of chords
in the figela processes $(S_{\alpha}(t))_{t\geq0}$.\vspace*{-1pt}

\begin{cor} \label{corofrag1} We have the following convergences:
\begin{longlist}[(ii)]
\item[(i)] If $\alpha=0$, $ e^{-t}\#S_{0}(t)
\overset{a.s.}{\underset{t\to\infty}{\longrightarrow}} \mathscr{E},$ where
$\mathscr{E}$ is exponentially distributed with parameter~1.
\item[(ii)] If $\alpha>0$,
$ t^{-1/\alpha}\#S_{\alpha}(t)\overset{a.s.}{\underset{t\to\infty
}{\longrightarrow}}
\frac{\Gamma(1/\alpha)}{\Gamma(2/\alpha)}$.\vspace*{-1pt}
\end{longlist}
\end{cor}

\begin{pf}(i) The case $p=0$ in assertion (ii) of Theorem \ref{outil}
gives the almost sure convergence
of the martingale $e^{-t}\#S_{0}(t)$. In fact, $ (\#S_{0}(t)
)_{t\geq0}$ is a Yule process of parameter $1$, which allows us to
identify the limit law (see \cite{AN72}, pages~127--130).

(ii) We first observe that
$\nu_{C}$ is conservative and thus $\mathscr{M}_\infty=1$ in the
notation of
Theorem \ref{outil}. The $\mathbb{L}^2$-convergence of $ t^{-1/\alpha
}\#(S_{\alpha}(t))$ toward a constant $K_{\nu_C}(\alpha,0)$
follows from Theorem \ref{outil}(iii) with $p=0$. From \cite{BD87},
Corollary~7, there is even almost sure convergence and the
constant $K_{\nu_C}(\alpha,0)$ is given by
$ K_{\nu_{C}}(\alpha,0) =\Gamma(1/\alpha)/\Gamma(2/\alpha)$.\vspace*{-1pt}
\end{pf}

A dissymmetry appears between the cases $\alpha=0$ and $\alpha>0$.
When $\alpha=0$, the number of chords grows exponentially with a
random multiplicative factor, but when $\alpha>0$ the number of chords
only grows like a power of $t$, with a deterministic multiplicative\vadjust{\goodbreak} factor.

\subsection{Fragments separating $1$ from a uniform point}
\label{fragsepuni}

Let $V$ be uniformly distributed over $\sun$ and independent of
$(S_{\alpha}(t),t\geq0,\alpha\geq0)$. Almost surely for every
$\alpha,t\geq0$, the points $1$ and $V$ do not belong to $\op
{Feet}(S_{\alpha}(t))$. Our goal is to establish a connection between
$H_{S_{\alpha}(t)}(1,V)$ [the height between $1$ and~$V$
in $S_{\alpha}(t)$] and a certain fragmentation process.

To this end, we first discuss the behavior of the figela process after
the appearance of
the first chord.
We briefly mentioned that the two fragments created by the first chord
of the figela process
can be viewed as two new disks by contracting the chord and that, after
the time of appearance
of the first chord, the process will
behave, independently in each of these two disks, as a rescaled copy of
the original process. Let us explain this in a more formal way. We fix
$\alpha\geq0$.

Let $[ab]$ be the first chord of the figela process
$(S_{\alpha}(t))_{t\geq0}$, which appears after an exponential time
$\tau$
with parameter $1$. We may write $a= e^{2i\pi U_{1}}, b= e^{2i\pi
U_{2}},$ where the pair $(U_{1},U_{2})$ has density $2\cdot \mathbf{1}_{\{0
<u_{1}<u_{2}<1\}}$ with respect to Lebesgue measure on $[0,1]^2$. Let
\[
M=1-(U_{2}-U_{1}),
\]
be the mass of the fragment of $S_{\alpha}(\tau)$ containing the
point $1$.

Define two mappings
$\psi_{U_{1},U_{2}}\dvtx[0,U_{1}]\cup[U_{2},1] \to[0,1]$ and $\phi
_{U_{1},U_{2}} \dvtx [U_{1}, U_{2}] \to[0,1]$ by setting
\begin{eqnarray}
\psi_{U_{1},U_{2}}(r) & =&
\cases{
\displaystyle
\frac{r}{M} ,&\quad if $0\leq r \leq
U_{1}$,\cr\displaystyle
\frac{r-(U_{2}-U_{1})}{M} ,&\quad if $U_{2}\leq
r\leq1$,
}
\nonumber\\
\phi_{U_{1},U_{2}}(r) &=&
\frac{r-U_{1}}{U_{2}-U_{1}} \qquad
\mbox{if } U_{1}\leq r \leq U_{2}.
\nonumber
\end{eqnarray}
Also let $ \Psi_{a,b}$ and $\Phi_{a,b}$ be the mappings corresponding
to $ \psi_{U_{1},U_{2}}$ and $ \phi_{U_{1},U_{2}}$
when $\sun$ is identified to $[0,1[$
\begin{eqnarray*}
\Psi_{a,b}(\exp(2i\pi r))&=& \exp(2i\pi\psi
_{U_{1},U_{2}}(r)), \qquad\mbox{if }r\in[0,U_1]\cup[U_2,1],\\
\Phi_{a,b}(\exp(2i\pi r))&=& \exp(2i\pi\phi
_{U_{1},U_{2}}(r)), \qquad\mbox{if }r\in[U_1,U_2].
\end{eqnarray*}

The first chord $[ab]$ creates two fragments. Let $R'$ the fragment (of
mass~$M$) containing $1$, and let $R''$ be the other fragment. For $t
\geq\tau$, we let $S^{(R')}_{\alpha}(t)$
[resp., $S^{(R'')}_{\alpha}(t)$] be the subset of $S_{\alpha
}(t)\setminus\{\{a,b\}\}$
consisting of all pairs $\{x,y\}$ such that the corresponding chord is
contained in $R'$ (resp., in $R''$).

\begin{lemma}
\label{keytool}
Let $\alpha\geq0$.
Conditionally on $(\tau, U_{1},U_{2})$, the pair of processes
\[
\bigl( \bigl(\Psi_{a,b} \bigl( S^{(R')}_{\alpha}(\tau+t) \bigr)
\bigr)_{t\geq0},
\bigl(\Phi_{a,b} \bigl( S^{(R'')}_{\alpha}(\tau+t) \bigr)
\bigr)_{t\geq0} \bigr)
\]
has the same distribution as
\[
\bigl( (S'_{\alpha}(M^\alpha t) )_{t\geq0},
\bigl(S''_{\alpha}\bigl((1-M)^\alpha t\bigr) \bigr)_{t\geq0} \bigr),
\]
where $S'_{\alpha}$ and $S''_\alpha$ are two independent copies of
the process $S_\alpha$.
\end{lemma}

This follows readily from our recursive construction of the figela process.

\begin{defi}
\label{sdc}
Let $S$ be a figela, and $x,y \in\sun\setminus\op{Feet}(S)$. We
call fragments \textit{separating} $x$ from $y$ in $S$, the fragments
of $S$ that intersect the chord~$[xy]$. These fragments are ranked
according to decreasing masses and denoted by
\[
R_{1}^{(x,y)}(S), R_{2}^{(x,y)}(S), R_{3}^{(x,y)}(S), \ldots.
\]
\end{defi}

\begin{figure}

\includegraphics{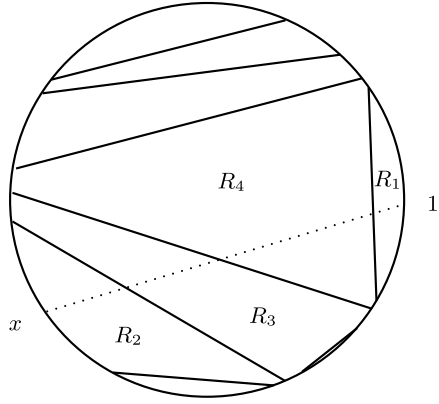}

\caption{A figela with 4 fragments $R_{1},R_{2},R_{3},R_{4}$
separating $1$ from $x$.}
\label{fig3}
\end{figure}

See Figure~\ref{fig3} for an example.

%

In order to state the main result of this subsection, we need one more
definition.
We let $\nu_D$ be the probability measure on $[0,1]^2$ defined by
\[
\int_{[0,1]^2}\nu_{D}(ds_{1},ds_{2}) F(s_{1},s_{2})= 2 \int
_{0}^1du u^2 F(u,0) + 4 \int_{1/2}^1 du u(1-u)F(u,1-u)
\]
for every nonnegative Borel function $F$.

The measure $\nu_{D}$ is interpreted as follows. Let $U,X_{1}$ and
$X_{2}$ be independent and uniformly distributed over $[0,1]$. The
point $U$ splits the interval $[0,1]$ in two parts, $[0,U[$ and
$]U,1]$. We keep each of these parts if and only if it contains at
least of of the two points $X_{1}$ or $X_{2}$. Then $\nu_D$
corresponds to the distribution of the masses of the remaining parts
ranked in decreasing order.

\begin{proposition} \label{uni} Let $V$ be a random variable uniformly
distributed over $\sun$ and independent of $(S_{\alpha}(t),t\geq0,
\alpha\geq0)$.\vadjust{\goodbreak} The sequence of masses of the fragments separating $1$
from $V$ in $S_{\alpha}(t)$, namely
\[
\mathcal{X}_{\alpha}(t) = \bigl(\mathrm{m}
\bigl(R_{1}^{(1,V)}(S_{\alpha}(t)) \bigr),\mathrm{m}
\bigl(R_{2}^{(1,V)}(S_{\alpha}(t)) \bigr),\ldots \bigr),
\]
is a fragmentation process with parameters $(\alpha, \nu_{D})$.
\end{proposition}

\begin{rek} Similarly as in Remark \ref{couplingremark}, the
coupling\vspace*{-1pt} of the processes $(S_{\alpha}(t),t\geq0)$ for $\alpha\geq
0$ induces the corresponding coupling of the processes $(\mathcal
{X}_{\alpha}(t),t\geq0)$ for $\alpha\geq0$.
\end{rek}

\begin{figure}

\includegraphics{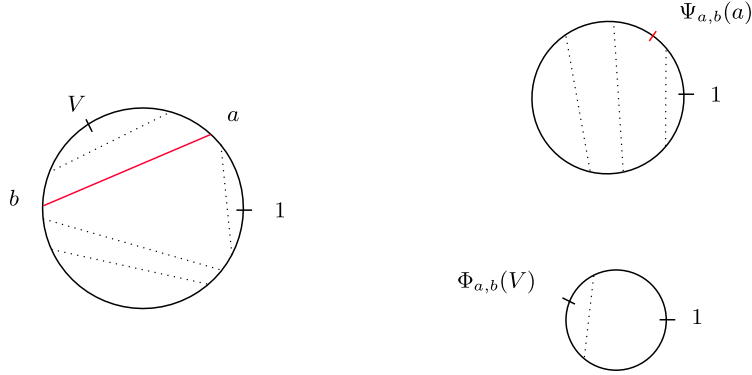}

\caption{Illustration of the proof in the case $V\in R''$.}
\label{fig4}
\end{figure}

\begin{pf*}{Proof of Proposition~\ref{uni}} We use the notation of the beginning of this subsection. Two
cases may occur.
\begin{longlist}[(2)]
\item[(1)] The point $V$ belongs to the fragment $R'$. Note that,
conditionally on the first chord $[ab]$ and on $\{ V\in R'\}$, $\Psi
_{a,b}(V)$ is uniformly distributed over $\sun$. Furthermore, the
future evolution of the process $\mathcal{X}_{\alpha}(t)$ after time~$\tau$
only depends on those chords that fall in the fragment $R'$
(and not on chords that fall in $R''$). More precisely, with the
notation of Lemma \ref{keytool},
the masses of the fragments of $S_\alpha(\tau+t)$ separating $1$ from
$V$ will be the same, up to the
mutiplicative factor $M$, as the
masses of the fragments of $\Psi_{a,b}(S^{(R')}_\alpha(\tau+t))$
separating $1$ from $\Psi_{a,b}(V)$.
By Lemma \ref{keytool}, conditionally on the event $\{ V\in R'\}$ and
on the pair $(\tau,M)$,
the process
$(\mathcal{X}_{\alpha}(\tau+t))_{t \geq0} $ has the same
distribution as
\[
(M \mathscr{X}_{\alpha}(M^\alpha t))_{t\geq0},
\]
where $(\mathscr{X}_{\alpha}(t))_{t\geq0}$ is a copy of $(\mathcal
{X}_{\alpha}(t))_{t\geq0}$,
which is independent of the pair $(\tau,M)$.
\item[(2)] The point $V$ belongs to the fragment $R''$ (see
Figure~\ref{fig4}). For
$t\geq\tau$, the fragments separating $V$ from $1$ in $S_{\alpha
}(t)$ will correspond either
to
fragments in the disk obtained from $R'$ by contracting the first chord
$[ab]$, provided these fragments separate $1$ from $\Psi_{a,b}(a)$, or to
fragments in the disk obtained
from $R''$ by contracting the first chord, provided these fragments
separate $\Phi_{a,b}(a)=1$ from $\Phi_{a,b}(V)$. An easy calculation
shows that, conditionally on $\{V\in R''\}$
and on $(\tau,M)$, the points $\Psi_{a,b}(a)$ and $\Phi_{a,b}(V)$
are independent and uniformly distributed over $\sun$.
Using Lemma \ref{keytool} once again, we get that the sequence of the
masses of
separating fragments contained in $R'$ at time $\tau+t$ has, as a
process in the variable $t$, the same distribution as $(M\mathscr
{X}_{\alpha}(M^\alpha t))_{t\geq0}$, where $\mathscr{X}_\alpha$
is an independent copy of $\mathcal{X}_{\alpha}$. A similar
observation holds for the
separating fragments in $R''$. Consequently, conditionally on the event
$\{V\in R''\}$
and on the pair $(\tau,M)$, the process
$ (\mathcal{X}_{\alpha}(\tau+t))_{t\geq0}$
has the same distribution as
\[
\bigl(M \mathscr{X}_{\alpha}(M^\alpha t)   \dotcupdisplay  (1-M)
\mathscr{X}'_{\alpha}\bigl((1-M)^\alpha t\bigr) \bigr)_{t\geq0},
\]
where $(\mathscr{X}_{\alpha}(t))_{t\geq0}$ and $(\mathscr
{X}'_{\alpha}(t))_{t\geq0}$ are independent copies of $(\mathcal
{X}_{\alpha}(t))_{t\geq0}$. Here the symbol $\dot\cup$ means that
we take the decreasing arrangement of the union of the two sequences.
\end{longlist}
Elementary calculations show that case $1$ occurs with probability
$2/3$ and that conditionally on this event the mass $M$ of the fragment
containing $1$ and~$V$ is distributed with density $3m^2$ on $[0,1]$.
Case $2$ occurs with probability $1/3$ and conditionally on that event
the mass of the largest fragment has density $12m(1-m)$ on $[1/2,1]$.
The preceding considerations then show that $(\mathcal{X}_{\alpha
}(t))_{t\geq0}$ is a fragmentation process with autosimilarity index
$\alpha$ and dislocation measure $\nu_{D}$ given as above.
\end{pf*}

%
%

In order to apply Theorem \ref{outil}
to the fragmentation process of Proposition~\ref{uni}, we must first
calculate the Malthusian exponent associated to $\nu_{D}$. From the
definition of $\nu_{D}$, we have for every $p\geq0$,
\[
\kappa_{\nu_{D}}(p) =1-2\int_{0}^1 du u^{p+2} -4 \int_{1/2}^1du
u(1-u) \bigl(u^p+(1-u)^p \bigr)
= \frac{p^2+3p-2}{p^2+5p+6}.
\]
Consequently, the only positive real $\beta^*$ such that $\kappa_{\nu
_{D}}(\beta^*)=0$ is
\[
\beta^* =\frac{\sqrt{17}-3}{2}.
\]
We also have $\kappa_{\nu_{D}}(0)=-1/3$.

Let $x \in\sun$. From now on, we will write
\[
\mathscr{M}^{(\alpha)}_{t}(x)= \sum_{i=1}^\infty\mathrm{m} \bigl(
R_{i}^{(1,x)}(S_{\alpha}(t)) \bigr)^{\beta^*}
\]
for the sum of the ${\beta^*}${th} powers of masses of the fragments
of $S_\alpha(t)$ separating $x$ from $1$. This makes sense since both
$1$ and $x$ a.s. do not belong to $\op{Feet}(S_\alpha(t))$.

By applying Theorem \ref{outil} to the fragmentation process $\mathcal
{X}_{\alpha}(t)$, we get:

\begin{cor} \label{coruni} Let $V$ be a random variable uniformly
distributed over $\sun$ and independent of $(S_{\alpha}(t),t\geq
0,\alpha\geq0)$. Then:\vadjust{\goodbreak}
\begin{longlist}[(iii)]
\item[(i)] The process $\mathscr{M}^{(\alpha)}_{t}(V)$ is a
uniformly integrable martingale and converges almost surely toward a
random variable $\mathscr{M}_{\infty}^V$ which does not depend on
$\alpha\geq0$.
Moreover $\mathscr{M}_\infty^V>0$ a.s., and $\mathbb{E} [ (\mathscr
{M}_\infty ^V)^q ]<\infty$ for every real $q\geq1$.
\item[(ii)] For every $\alpha>0$, there exists a constant $K_{\nu
_{D}}(\alpha)$ such that
\[
t^{-\beta^*/\alpha} H_{S_{\alpha}(t)}(1,V) \underset{t\to\infty
}{\overset{\mathbb{L}^2}{\longrightarrow}} K_{\nu_{D}}(\alpha)
\mathscr{M}_{\infty}^V.
\]
\item[(iii)] There exists a positive random variable $\mathscr
{H}_{0}^V$ such that
\[
e^{-t/3} H_{S_0(t)}(1,V)\underset{t\to\infty}{\overset
{a.s.}{\longrightarrow}} \mathscr{H}_{0}^V.
\]
More generally, for every $p\geq0$, there exists
a positive random variable $\mathscr{H}_{p}^V$ such that
\[
e^{t\kappa_{\nu_{D}}(p)} \sum_{i=0}^\infty\mathrm{m}
\bigl(R_{i}^{(1,V)}(S_{0}(t)) \bigr)^p
\underset{t\to\infty}{\overset{a.s.}{\longrightarrow}} \mathscr{H}_{p}^V.
\]
\end{longlist}
\end{cor}

\begin{rek}
\label{technicalrem}
The convergence in (ii) can be reinforced in the following way. For
every $\delta\in\,]0,1[$,
%
%
\begin{equation}
\label{technicalpoint}
\lim_{t\to\infty} \mathbb{E} \Bigl[ \sup_{\delta t\leq s\leq t}
\bigl|s^{-\beta ^*/\alpha} H_{S_\alpha(s)}(1,V) - K_{\nu_D}(\alpha)
\mathscr {M}_\infty^V \bigr|^2 \Bigr]=0.
\end{equation}
To see this, fix $\varepsilon\in\,]0,1[$ and choose a subdivision
$\delta=\delta_0<\delta_1<\cdots<\delta_k=1$ of $[\delta,1]$ such that
$(\delta_{i+1}/\delta_i)^{\beta^*/\alpha}<1+\varepsilon$ for every
$0\leq i\leq k-1$. Since the function
$s\mapsto H_{S_\alpha}(s)(1,V) $ is nondecreasing, we have
\begin{eqnarray*}
&&\sup_{\delta t\leq s\leq t} \bigl(s^{-\beta^*/\alpha} H_{S_\alpha
(s)}(1,V) - K_{\nu_D}(\alpha) \mathscr{M}_\infty^V \bigr)\\
&& \qquad \leq\sup_{0\leq i\leq k-1} \bigl( (\delta_it)^{-\beta
^*/\alpha} H_{S_\alpha(\delta_{i+1}t)}(1,V) - K_{\nu_D}(\alpha)
\mathscr{M}_\infty^V \bigr)\\
&& \qquad \leq(1+\varepsilon)
\sup_{0\leq i\leq k-1} \bigl( (\delta_{i+1}t)^{-\beta^*/\alpha}
H_{S_\alpha(\delta_{i+1}t)}(1,V) - K_{\nu_D}(\alpha) \mathscr
{M}_\infty^V \bigr)\\
&& \qquad  \quad {}
+ \varepsilon K_{\nu_D}(\alpha) \mathscr{M}_\infty^V.
\end{eqnarray*}
Similar manipulations give
\begin{eqnarray*}
&&\sup_{\delta t\leq s\leq t} \bigl(K_{\nu_D}(\alpha) \mathscr
{M}_\infty^V - s^{-\beta^*/\alpha} H_{S_\alpha(s)}(1,V) \bigr)\\
&& \qquad \leq\sup_{0\leq i\leq k-1} \bigl(K_{\nu_D}(\alpha) \mathscr
{M}_\infty^V -
(\delta_{i+1}t)^{-\beta^*/\alpha} H_{S_\alpha(\delta
_{i}t)}(1,V) \bigr) \\
&& \qquad \leq(1-\varepsilon)
\sup_{0\leq i\leq k-1} \bigl(K_{\nu_D}(\alpha) \mathscr{M}_\infty
^V - (\delta_{i}t)^{-\beta^*/\alpha} H_{S_\alpha(\delta
_{i}t)}(1,V) \bigr)\\
&& \qquad  \quad {}
+ \varepsilon K_{\nu_D}(\alpha) \mathscr{M}_\infty^V.
\end{eqnarray*}
It follows that
\begin{eqnarray*}
&&\sup_{\delta t\leq s\leq t} \bigl|s^{-\beta^*/\alpha} H_{S_\alpha
(s)}(1,V) - K_{\nu_D}(\alpha) \mathscr{M}_\infty^V \bigr|\\
&& \qquad \leq2 \sup_{0\leq i\leq k} \bigl|(\delta_{i}t)^{-\beta
^*/\alpha} H_{S_\alpha(\delta_{i}t)}(1,V) -
K_{\nu_D}(\alpha) \mathscr{M}_\infty^V \bigr| + \varepsilon
K_{\nu_D}(\alpha) \mathscr{M}_\infty^V.
\end{eqnarray*}
Using property (ii) in Corollary \ref{coruni}, we now get
\[
\limsup_{t\to\infty} \mathbb{E}\Bigl [ \sup_{\delta t\leq s\leq t}
\bigl|s^{-\beta^*/\alpha} H_{S_\alpha(s)}(1,V)\,{-}\,K_{\nu_D}(\alpha)
\mathscr{M}_\infty^V \bigr|^2 \Bigr]\,{\leq}\,2 \varepsilon^2 K_{\nu_D}(\alpha)^2 \mathbb{E} [ (\mathscr
{M}_\infty^V)^2 ],
\]
and (\ref{technicalpoint}) follows since $\varepsilon$ was arbitrary.
\end{rek}

\subsection{Fragments separating $1$ from a deterministic point}
\label{fragdeter}
We now aim at an analogue of the last corollary when $V$ is replaced by
a deterministic point~$x$ in $\sun$.
We will use the position of the first chord to provide the randomness
that we need to reduce the proof to the statement of Corollary~\ref
{coruni}. We start by computing explicitly the distributions of certain
quantities that arise when describing the evolution of the process
after the creation of the first chord.
We use the notation of the beginning of the previous subsection.

We fix $r\in\,]0,1[$ and write $x=e^{2i\pi r}$. Consider first the case when
$x\in R''$, or equivalently
$U_{1}<r<U_{2}$. We then set $Y_{1}= \psi_{U_{1},U_{2}}(U_{1}) = \frac
{U_{1}}{M}$, which represents the position of the distinguished point,
corresponding to the endpoints of the first chord, in the disk obtained
from $R'$ by contracting the first chord. Similarly, $Y_{2}=\phi
_{U_{1},U_{2}}(r)=\frac{r-U_{1}}{1-M}$ gives the position of the
distinguished point corresponding to $x$ in the the disk obtained from
$R''$ by the same contraction.

In the case when $r\in \,]0,U_{1}[\,\cup\,]U_{2},1[$ (or equivalently
$x\in R'$), we take $Y_{2}=0$ and we let $Y_{1}= \psi
_{U_{1},U_{2}}(r)$ be the position of the point corresponding to $x$ in
the new disk obtained from $R'$ by contracting the first chord.

We first evaluate the density of the pair $(Y_{1},Y_{2})$
on the event $\{x\in R''\}=\allowbreak\{U_1\,{<}\,r\,{<}\,U_2\}$. We have, for any nonnegative
measurable function $f$ on~$[0,1]^2$,
\begin{eqnarray*}
&&\mathbb{E} \bigl[ f(Y_{1},Y_{2}) {\mathbf1}_{\{U_1<r<U_2\}} \bigr]\\
&& \qquad =2\iint_{[0,1]^2} du_{1}\,du_{2} \mathbf{1}_{\{
u_{1}<r<u_{2}\}}
f \biggl(\frac{u_1}{1-(u_2-u_1)},\frac{r-u_1}{u_2-u_1} \biggr)\\
&& \qquad =2r(1-r) \int_{0}^1 \!\! \int_{0}^1 ds_{1}\,ds_{2}\\
&& \qquad  \quad
{}\times f
\biggl(\frac{rs_{1}}{rs_{1}+(1-r)(1-s_{2})},\frac
{r(1-s_{1})}{r(1-s_{1})+(1-r)s_{2}} \biggr).
\end{eqnarray*}
From the obvious change of variables and after tedious calculations, we get
%
%
\begin{equation}
\label{eq1}
\mathbb{E} \bigl[ f(Y_{1},Y_{2}) {\mathbf1}_{\{U_1<r<U_2\}} \bigr]
=2\iint_{\mathscr{D}_{r}} dr_{1}\,dr_{2} \frac
{|r_{1}-r||r_{2}-r|}{|r_{1}-r_{2}|^3} f(r_{1},r_{2}),
\end{equation}
where $\mathscr{D}_{r}$ is the set $\mathscr{D}_{r}=
([r,1]\times[0,r] ) \cup([0,r]\times[r,1] ) $.
Also note that, on the event $\{U_1<r<U_2\}$, we have $MY_{1}+(1-M)Y_{2}=r$,
and thus $M=\frac{r-Y_2}{Y_1-Y_2}$.\vspace*{1pt}

We can similarly compute the distribution of $Y_{1}$ on the event \mbox{$\{
U_1\,{<}\,U_2\,{<}\,r\}$}. For any nonnegative measurable function $f$ on $[0,1]$,
%
%
\begin{eqnarray}
\label{eq2}
\mathbb{E} \bigl[ f(Y_{1}) {\mathbf1}_{\{U_1<U_2<r\}} \bigr]
& =& 2\iint_{[0,1]^2} du_{1}\,du_{2} \mathbf{1}_{\{u_{1}<u_{2}<r\}}
f \biggl(\frac{r-(u_2-u_1)}{1-(u_2-u_1)} \biggr)\notag\\
& = & 2r^2 \int_{0}^1\int_{0}^1 ds_{1}\,ds_{2} \mathbf{1}_{\{
s_{1}<s_{2}\}} f \biggl(\frac
{r(1-(s_{2}-s_{1}))}{1-r(s_{2}-s_{1})} \biggr)\\
& =& 2(1-r)^2 \int_{0}^r dr_{1}\frac{r_{1}}{(1-r_{1})^3} f(r_1).
\notag
\end{eqnarray}
Also notice that $M=\frac{1-r}{1-Y_{1}}$ on the event $\{U_1<U_2<r\}$.

A similar calculation, or a symmetry argument, shows that the
distribution of $Y_{1}$ on the event $\{r<U_1<U_2\}$ is given by
%
%
\begin{equation}
\label{eq3}
\mathbb{E} \bigl[ f(Y_{1}) {\mathbf1}_{\{r<U_1<U_2\}} \bigr] = 2r^2 \int_{r}^1
dr_{1}\frac{1-r_{1}}{r_{1}^3} f(r_1).
\end{equation}
Note that $M=\frac{r}{Y_{1}}$ on $\{r<U_1<U_2\}$.

We can now state and prove the main result of this section.

\begin{theo}\label{tech} Let $x\in\sun\setminus\{1\}$.
\begin{longlist}[(iii)]
\item[(i)] For every $\alpha\geq0$, the process $\mathscr
{M}_{t}^{(\alpha)}(x)$ converges almost surely toward a random
variable $\mathscr{M}_{\infty}(x)$ which does not depend on $\alpha$.
\item[(ii)] We have $\mathscr{M}_{\infty}(x)>0$ a.s. and $\mathbb
{E} [ \mathscr{M}_\infty(x)^q ]<\infty$ for every $q\geq1$.
\item[(iii)] For every $\alpha>0$, we have
\[
t^{-\beta^*/\alpha}H_{S_{\alpha}(t)}(1,x) \overset{(\mathbb
{P})}{\underset{t\to\infty}{\longrightarrow}} K_{\nu_{D}}(\alpha)
\mathscr{M}_{\infty}(x),
\]
where the constant $K_{\nu_{D}}(\alpha)$ is the same as in Corollary
\ref{coruni}.
\item[(iv)] There exists a positive random variable $\mathscr
{H}_{0}(x)$ such that
\[
e^{-t/3} H_{S_0(t)}(1,x)\underset{t\to\infty}{\overset
{a.s.}{\longrightarrow}} \mathscr{H}_{0}(x).
\]
More generally, for every $p\geq0$, there exists a positive random
variable~$\mathscr{H}_{p}(x)$ such that
\[
e^{t\kappa_{\nu_{D}}(p)} \sum_{i=0}^\infty\mathrm{m}
\bigl(R_{i}^{(1,x)}(S_{0}(t)) \bigr)^p \overset{a.s.}{\underset{t\to
\infty}{\longrightarrow}} \mathscr{H}_{p}(x).
\]
\end{longlist}
\end{theo}

\begin{pf} As previously, we write $x=e^{2i\pi r}$, where $r\in\,]0,1[$.
To simplify notation, we also set, for every $t\geq0$, and every
$\alpha\geq0$,
\[
\mathcal{X}^x_{\alpha}(t) = \bigl(\mathrm{m}
\bigl(R_{1}^{(1,x)}(S_{\alpha}(t)) \bigr),\mathrm{m}
\bigl(R_{2}^{(1,x)}(S_{\alpha}(t)) \bigr),\ldots \bigr).
\]

Fix $\alpha\geq0$. Consider first the case when $x$ belongs to $R'$.
After time $\tau$, the fragments separating $1$ from $x$ will
correspond to fragments separating $1$ from $\Psi_{a,b}(x)$ in the
disk obtained from $R'$ by contracting the first chord. If $F$ is a
nonnegative measurable function on the Skorokhod space
${\mathbb D}([0,\infty[,{\mathcal S}^\downarrow)$, Lemma \ref
{keytool} gives
%
%
\begin{equation}
\label{techeq1} \quad
\mathbb{E} \bigl[ F \bigl(\bigl({\mathcal X}^x_{\alpha}(\tau+t)\bigr)_{t\geq0} \bigr)
{\mathbf1}_{\{x\in R'\}} \bigr]
=\mathbb{E} \bigl[ F \bigl(\bigl(M\widetilde{\mathcal X}^{\Psi_{a,b}(x)}_{\alpha
}(M^\alpha t)\bigr)_{t\geq0} \bigr){\mathbf1}_{\{x\in R'\}} \bigr],
\end{equation}
where, for every $y\in\sun\setminus\{1\}$, the process
$(\widetilde{\mathcal X}^{y}_{\alpha}(t))_{t\geq0}$ is defined from
an independent copy
$(\widetilde S_\alpha(t))_{t\geq0}$ of $(S_\alpha(t))_{t\geq0}$, in
the same way as
$({\mathcal X}^{y}_{\alpha}(t))_{t\geq0}$ is defined from $(S_\alpha
(t))_{t\geq0}$. Note that
$\Psi_{a,b}(x)=\exp(2i\pi Y_1)$ in the notation introduced before the theorem.
From formulas (\ref{eq2}) and (\ref{eq3}) and the relations between
$M$ and $Y_1$, we get
\begin{eqnarray*}
&&\mathbb{E} \bigl[ F \bigl(\bigl({\mathcal X}^x_{\alpha}(\tau+t)\bigr)_{t\geq0} \bigr)
{\mathbf1}_{\{x\in R'\}} \bigr]\\
&& \qquad =2(1-r)^2 \int_0^r dr_1 \frac{r_1}{(1-r_1)^3}\\
&& \qquad   \quad
{}\times\mathbb{E}
\biggl[ F \biggl( \biggl( \biggl(\frac{1-r}{1-r_1} \biggr) \widetilde{\mathcal X}^{\exp(2i\pi
r_1)}_{\alpha} \biggl( \biggl(\frac {1-r}{1-r_1} \biggr)^\alpha t \biggr) \biggr)_{t\geq0} \biggr) \biggr]\\
&& \qquad \quad{} + 2r^2 \int_r^1 dr_1 \frac{1-r_1}{r_1^3} \mathbb
{E} \biggl[ F \biggl( \biggl( \biggl(\frac{r}{r_1} \biggr) \widetilde{\mathcal X}^{\exp(2i\pi
r_1)}_{\alpha} \biggl( \biggl(\frac {r}{r_1} \biggr)^\alpha t \biggr) \biggr)_{t\geq0} \biggr) \biggr].
\end{eqnarray*}
Let $U$ be uniformly distributed over $[0,1]$ and independent of
$(\widetilde S_\alpha(t))_{t\geq0}$.
By the preceding display, the conditional distribution of $({\mathcal
X}^x_{\alpha}(\tau+t))_{t\geq0}$
given that $x\in R'$ is absolutely continuous (even with a bounded
density) with respect to that of the process
\begin{eqnarray*}
&&\biggl({\mathbf1}_{\{U<r\}} \frac{1-r}{1-U}
\widetilde{\mathcal X}^{\exp(2i\pi U)}_{\alpha} \biggl( \biggl(\frac
{1-r}{1-U} \biggr)^\alpha t \biggr)\\
&& \qquad \quad
+ {\mathbf1}_{\{U>r\}} \frac{r}{U}
\widetilde{\mathcal X}^{\exp(2i\pi U)}_{\alpha} \biggl( \biggl(\frac
{r}{U} \biggr)^\alpha t \biggr)
\biggr)_{t\geq0}.
\end{eqnarray*}
Since $V=\exp(2i\pi U)$ is uniformly distributed on $\sun$ and
independent of $(\widetilde S_\alpha(t))_{t\geq0}$, we can apply
Corollary \ref{coruni} to get asymptotics for the process
in the last display. It follows that the almost sure convergences in
parts~(i) and (iv) of the proposition hold on the event $\{x\in R'\}$.
Moreover the variable~$\mathscr{M}_{\infty}(x)$ obtained as the
almost sure limit of $\mathscr{M}_{t}^{(\alpha)}(x)$
(only on the event $\{x\in R'\}$ for the moment) does not depend on the
choice of $\alpha\geq0$.
To see this, note that if we fix two values $\alpha\geq0$ and $\alpha
'\geq0$, the
preceding absolute continuity property holds in a similar form for the pair
$(({\mathcal X}^x_{\alpha}(\tau+t))_{t\geq0},({\mathcal X}^x_{\alpha
'}(\tau+t))_{t\geq0})$.
Then it suffices to use the fact that the limiting variable ${\mathscr
M}_\infty^V$ in
Corollary~\ref{coruni}(i) does not depend on the choice of $\alpha
\geq0$.\looseness=-1

The justification of property (iii) of the theorem (still on the event
$\{x\in R'\}$) is a bit trickier because we do not have almost sure
convergence in
Corollary \ref{coruni}(ii). We need the reinforced version of
Corollary \ref{coruni}(ii) provided by Remark
\ref{technicalrem}. We observe that, if $U>r$, the quantity
$(r/U)^\alpha$ is bounded above by $1$
and bounded below by $r^\alpha$, so that (\ref{technicalpoint}) gives
\[
t^{-\beta^*/\alpha}
H_{\widetilde S_\alpha((r/U)^\alpha t)}(1,\exp(2i\pi U))
\overset{(\mathbb{L}^2)}{\underset{t\to\infty}{\longrightarrow}}
\biggl(\frac{r}{U} \biggr)^{\beta^*} K_{\nu_D}(\alpha) \tilde
{\mathscr M}_\infty^V
\]
on the event $\{U>r\}$ (with an obvious notation for $ \tilde{\mathscr
M}_\infty^V$). A similar observation holds for
the asymptotics of $H_{\widetilde S_\alpha(((1-r)/(1-U))^\alpha
t)}(1,\exp(2i\pi U))$ on the event $\{U<r\}$. By combining
both asymptotics and using the absolute continuity relation mentioned
above, we get that the convergence in
probability in assertion (iii) of the proposition holds on the event $\{
x\in R'\}$.

Let us turn to the case where $x$ belongs to $R''$. From Lemma \ref
{keytool}, we have
%
%
\begin{eqnarray}
\label{techeq2}
&& \quad \mathbb{E} \bigl[ F \bigl(\bigl({\mathcal X}^x_{\alpha}(\tau+t)\bigr)_{t\geq0} \bigr)
{\mathbf1}_{\{x\in R''\}} \bigr]\nonumber
\\[-1pt]
&& \quad  \qquad =\mathbb{E} \bigl[ F \bigl( \bigl(M\widetilde{\mathcal X}^{\Psi
_{a,b}(a)}_{\alpha }(M^\alpha t) \\[-1pt]
&& \quad  \qquad\hphantom{=\mathbb{E} \bigl[ F \bigl( (} {}\dotcupdisplay (1-M) {\bar{\mathcal
X}}_\alpha^{\Phi _{a,b}(x)}\bigl((1-M)^\alpha t\bigr) \bigr)_{t\geq0} \bigr) {\mathbf
1}_{\{x\in R''\}} \bigr],
\nonumber
\end{eqnarray}
where $\widetilde X^y_\alpha$ and $\bar X^y_\alpha$ are defined in
terms of
two independent copies $\widetilde S_\alpha$ and~$\bar S_\alpha$ of
$S_\alpha$
(and the notation $\dot\cup$ has the same meaning as in the proof of
Proposition \ref{uni}).

Using now formula (\ref{eq1}), we obtain
\begin{eqnarray*}
&&\hspace*{-4pt}\mathbb{E} \bigl[ F \bigl(\bigl({\mathcal X}^x_{\alpha}(\tau+t)\bigr)_{t\geq0} \bigr)
{\mathbf1}_{\{x\in R''\}} \bigr]\notag\\[-1pt]
&&\hspace*{-4pt} \qquad = 2
\int\!\!\int_{\mathscr{D}_r} dr_1\,dr_2\frac{|r_1-r| |r_2-r|}{|r_1-r_2|^3}
\\[-1pt]
&&\hspace*{-4pt}
\qquad  \quad  {} \times\mathbb{E} \biggl[ F \biggl( \biggl( \biggl( \frac{r-r_2}{r_1-r_2}
\biggr)\widetilde{\mathcal X}_\alpha^{\exp(2i\pi r_1)}  \biggl( \biggl( \frac
{r-r_2}{r_1-r_2} \biggr)^\alpha t \biggr)\\[-1pt]
&& \hphantom{\times\mathbb{E} \biggl[ F \biggl( \biggl(}
 \qquad  \quad   \dotcupdisplay \biggl(\frac{r_1-r}{r_1-r_2} \biggr) \bar
{\mathcal X}_\alpha^{\exp(2i\pi r_2)} \biggl( \biggl( \frac{r_1-r}{r_1-r_2}
\biggr)^\alpha t \biggr) \biggr)_{t\geq0} \biggr) \biggr].\nonumber
\end{eqnarray*}
Hence, if $U$ and $U'$ are two independent variables uniformly
distributed over $[0,1]$
and independent of $(\widetilde S_\alpha,\bar S_\alpha)$, the
distribution of $({\mathcal X}^x_{\alpha}(\tau+t))_{t\geq0}$ knowing
that $x\in R''$
is absolutely continuous with respect to the distribution of
\begin{eqnarray*}
&&\biggl( \biggl( \frac{r-U'}{U-U'} \biggr)\widetilde{\mathcal X}_\alpha
^{\exp(2i\pi U)}
\biggl( \biggl( \frac{r-U'}{U-U'} \biggr)^\alpha t \biggr)\\[-1pt]
&& \qquad
\dotcupdisplay \biggl(\frac{U-r}{U-U'} \biggr) \bar{\mathcal X}_\alpha
^{\exp(2i\pi U')}
\biggl( \biggl( \frac{U-r}{U-U'} \biggr)^\alpha t \biggr) \biggr)_{t\geq0}
\end{eqnarray*}
conditionally on $(U-r)(U'-r)<0$. As in the case $x\in R'$, we see that
the almost sure convergences in assertions (i) and (iv) of the
proposition,\vadjust{\goodbreak} on the
event $\{x\in R''\}$, follow from the analogous convergences in
Corollary \ref{coruni}. By the same argument as in
the case $x\in R'$, the almost sure limit ${\mathscr M}_\infty(x)$ in
(i) does not depend on the choice
of $\alpha\geq0$.

To get the convergence in probability in assertion (iii), we again use
Remark~\ref{technicalrem}. The point is that the quantities
$((r-U')/(U-U'))^\alpha$ and $((U-r)/(U-U'))^\alpha$, which are
bounded above by $1$ [recall that we condition on $(U-r)(U'-r)<0$],
are also bounded below by $\delta>0$ except on a~set of small
probability. As in the case $x\in R'$, the desired result follows
from~(\ref{technicalpoint}).

It remains to prove (ii). The property $\mathscr{M}_\infty(x)>0$ a.s.
is immediate from the analogous property in
Corollary \ref{coruni} and our absolute continuity argument. Then, by
applying formulas (\ref{techeq1}) and (\ref{techeq2})
with a suitable choice of the function $F$, we get, for every
nonnegative measurable function
$f$ on ${\mathbb R}_+$,
\begin{eqnarray*}
&&\mathbb{E} [ f(\mathscr{M}_\infty(x)) ]\\
&& \qquad= \mathbb{E} \bigl[ f(\mathscr{M}_\infty (x)){\mathbf1}_{\{
x\in R'\}} \bigr]
+ \mathbb{E} \bigl[ f(\mathscr{M}_\infty(x)){\mathbf1}_{\{x\in R''\}} \bigr]\\
&& \qquad= \mathbb{E} \bigl[ f(M^{\beta^*}\tilde{\mathscr M}_\infty
(e^{2i\pi Y_1})){\mathbf1}_{\{x\in R'\}} \bigr]\\
&& \qquad  \quad {}
+ \mathbb{E} \bigl[ f\bigl(M^{\beta^*}\tilde{\mathscr M}_\infty(e^{2i\pi
Y_1})+(1-M)^{\beta^*}\bar{\mathscr M}_\infty(e^{2i\pi Y_2})\bigr)
{\mathbf1}_{\{x\in R''\}} \bigr],
\end{eqnarray*}
where $\tilde{\mathscr M}_\infty$ and $\bar{\mathscr M}_\infty$ are
the obvious analogues of $\mathscr{M}_\infty$
when $S_\alpha$ is replaced by $\widetilde S_\alpha$ and $\bar
S_\alpha$, respectively. Set $\mathscr{M}_\infty(1)=0$. We
have obtained the identity in distribution
%
%
\begin{equation} \label{idloi} \mathscr{M}_{\infty}(x)\overset{(d)}{=}
M^{\beta^*}\mathscr{M}'_{\infty}(e^{2i\pi Y_{1}}) +(1-M)^{\beta
^*}\mathscr{M}''_{\infty}(e^{2i\pi Y_{2}}),
\end{equation}
where $\mathscr{M}'_{\infty}$ and $\mathscr{M}_{\infty}''$ are two
independent copies of $\mathscr{M}_{\infty}$, and the pair $(\mathscr
{M}'_{\infty},\allowbreak\mathscr{M}''_{\infty})$ is also independent of
$(Y_1,Y_2)$. However, from the explicit formulas (\ref{eq1}), (\ref
{eq2}) and (\ref{eq3}), it is easy to
verify that both the density of the law of~$Y_1$ and the density of the
law of $Y_2$ conditional on $\{Y_2\not= 0\}$ are
bounded above by a~constant depending on $x$ [even though the joint
density of the pair $(Y_1,Y_2)$ conditional
on $\{Y_2\not= 0\}$ is unbounded]. By Corollary \ref{coruni} we know
that, if~$U$ is uniformly
distributed over $[0,1]$ and independent of the figela process, we have
$\mathbb{E} [ \mathscr{M}_\infty(e^{2i\pi U})^q ]<\infty$
for every $q\geq1$. The analogous property for~$\mathscr{M}_{\infty
}(x)$ then follows from
(\ref{idloi}) and the preceding observations.
\end{pf}

\begin{rek}
By rotational invariance of the model, point $1$ can be replaced by any
point of $\sun$ in
Theorem \ref{tech}.
\end{rek}

\section{Estimates for moments and the continuity of the height
process}\label{sec4}

\subsection{Estimates for moments}

We first state a proposition
giving estimates for the moments of the increments of the process
$\mathscr{M}_{\infty}(x)$.
These estimates will allow us to apply Kolmogorov's continuity criterion
in order to get information on the H\"{o}lder continuity properties of
this process.
Recall that we take $\mathscr{M}_{\infty}(1)=0$ by\vadjust{\goodbreak} convention.

\begin{proposition}\label{moments} For every $\varepsilon>0$ and
every integer $p\geq1$, there exists a constant $M_{\varepsilon
,p}\geq0$ such that, for every $u\in[0,1]$ we have
\[
\mathbb{E} [ \mathscr{M}_{\infty}(e^{2i\pi u})^p ] \leq
M_{\varepsilon
,p}\bigl(u(1-u)\bigr)^{p\beta^*-\varepsilon}.
\]
In the special case $p=1$, we have
\[
\mathbb{E} [ \mathscr{M}_{\infty}(e^{2i\pi u}) ] = \frac{\Gamma
(2+2\beta
^*)}{\Gamma(1+\beta^*)^2}\bigl(u(1-u)\bigr)^{\beta^*}.
\]
\end{proposition}

The proof of the proposition is given in the next two subsections. This
proof relies
on the identity in distribution (\ref{idloi}) derived in the preceding
proof. Using this identity and
formulas (\ref{eq1}), (\ref{eq2}) and (\ref{eq3}), we will obtain
integral equations for the moments of $(\mathscr{M}_{\infty}(x),x\in
\sun)$. We can explicitly solve the integral equation corresponding to
the first moment.
We then use Gronwall's lemma to investigate the behavior of higher
moments when $x\in\sun$ is close to $1$.

For every integer $p\geq1$ and every $r\in[0,1]$, we set
\[
m_{p}(r) =\mathbb{E} [ \mathscr{M}_{\infty}(e^{2i\pi r})^p ].
\]

\subsection{The case $p=1$}
Let $r\in\,]0,1[$.
Thanks to the identity in distribution~(\ref{idloi}) and to formulas
(\ref{eq1}), (\ref{eq2}) and (\ref{eq3}), we obtain the integral equation
%
%
\begin{eqnarray}
\label{sump=1}
m_{1}(r) &= & 2(1-r)^2\int_{0}^rdr_{1} \frac{r_{1}}{(1-r_{1})^3}
\biggl(\frac{1-r}{1-r_{1}} \biggr)^{\beta^*}m_{1}(r_{1})\nonumber\\[-1pt]
&&{}
+2r^2\int_{r}^{1}dr_{1} \frac{1-r_{1}}{r_{1}^3} \biggl(\frac
{r}{r_{1}} \biggr)^{\beta^*}m_{1}(r_{1})\nonumber
\\[-10pt]
\\[-10pt]
&&{}+ 2\iint_{\mathscr{D}_{r}} dr_{1}\,dr_{2} \frac
{|r_{1}-r||r_{2}-r|}{|r_{1}-r_{2}|^3}\nonumber\\[-1pt]
&&
\quad{}\times \biggl( \biggl(\frac
{r-r_{2}}{r_{1}-r_{2}} \biggr)^{\beta^*}m_{1}(r_{1}) + \biggl(\frac
{r_{1}-r}{r_{1}-r_{2}} \biggr)^{\beta^*}m_{1}(r_{2}) \biggr).
\nonumber
\end{eqnarray}
We can rewrite the first two terms in the sum of the right-hand side in
the form
\[
2\!\int_0^{r}\!dr_1\!\biggl(\!\frac{1}{1\,{-}\,r_1}\,{-}\,1 \!\biggr)\!\biggl(\!\frac
{1\,{-}\,r}{1\,{-}\,r_1}\!\biggr) ^{\beta^*+2}
m_1(r_1)\,{+}\,2\!\int_r^{1}\!dr_1\!\biggl(\!\frac{1}{r_1}\,{-}\,1\!\biggr)\biggl (\!\frac
{r}{r_1} \!\biggr) ^{\beta^*+2} m_1(r_1).
\]
As for the third term, we observe that
\begin{eqnarray*}
&& \int_{0}^r dr_{1}\int_{r}^1 dr_{2} \frac
{(r-r_{1})(r_{2}-r)}{(r_{2}-r_{1})^3} \biggl(\frac
{r_2-r}{r_{2}-r_{1}} \biggr)^{\beta^*}m_{1}(r_{1}) \\[-1pt]
&& \qquad = \int_{0}^r dr_{1} m_{1}(r_{1}) (r-r_{1}) \int_{r}^1
dr_{2}
\biggl(\frac{1}{r_{2}-r_1} \biggr)^2\biggl (\frac
{r_2-r}{r_{2}-r_{1}} \biggr)^{\beta^*+1}\\[-1pt]
&& \qquad = \int_{0}^r dr_{1} m_{1}(r_{1}) \frac{1}{\beta^*+2}
\biggl( \frac{1-r}{1-r_{1}} \biggr) ^{\beta^*+2} ,
\end{eqnarray*}
where we made the change of variables $u=\frac{r_2-r}{r_2-r_1}$ to
compute the integral in~$dr_2$. It follows that the third term in the
right-hand side of (\ref{sump=1}) is equal to
\[
\frac{4}{\beta^*+2} \biggl(\int_{0}^r dr_{1} \biggl( \frac
{1-r}{1-r_{1}} \biggr) ^{\beta^*+2} m_{1}(r_{1})
+ \int_{r}^1 dr_{1} \biggl( \frac{r}{r_{1}} \biggr) ^{\beta^*+2}
m_{1}(r_{1}) \biggr) .
\]
Summarizing, we obtain that the function $(m_1(r),r\in\,]0,1[)$ solves
the integral equation
%
%
\begin{equation}\label{eqm1}
m_{1}(r)=\int_{0}^1 du g_{r}(u) m_{1}(u),
\end{equation}
where, for every $r\in\,]0,1[$,
\begin{eqnarray*}
g_{r}(u)& =& \mathbf{1}_{\{0< u <r\}} \biggl(\frac{1-r}{1-u}
\biggr)^{2+\beta^*} \biggl(\frac{2}{1-u}-\frac{2\beta^*}{\beta^*+2}
\biggr)\\
&&{}
+ \mathbf{1}_{\{r \leq u < 1\}}\biggl (\frac{r}{u} \biggr)^{2+\beta
^*}\biggl (\frac{2}{u}-\frac{2\beta^*}{\beta^*+2} \biggr),
\end{eqnarray*}
is a positive function on $]0,1[$. Elementary calculations, using the
fact that $ ( \beta^* )^2+3\beta^* -2=0$, show that $\int
_{0}^1 g_{r}(u)\,dr =1$, for every $u\in\,]0,1[$.

Let $N$ be the operator that maps a function $f \in\mathbb
{L}^1(]0,1[,dr)$ to the function
\[
N(f)(r)=\int_{0}^1 du g_{r}(u)f(u).
\]
Then $N$ is a contraction: If $f_{1},f_{2}\in\mathbb{L}^1(]0,1[,dr)$,
we have
\begin{eqnarray*}
\int_{0}^1dr |N(f_{1})(r)-N(f_{2})(r)|
&\leq&\int_{0}^1dr \int_{0}^1 du g_{r}(u)|f_{1}(u)-f_{2}(u)|\\
&=& \int_{0}^1du |f_{1}(u)-f_{2}(u)|.
\end{eqnarray*}
The first inequality in the last display is strict unless $f_{1}-f_{2}$
has a.e. a~constant sign.
It follows that there can be at most one nonnegative function $f \in
\mathbb{L}^1(]0,1[,dr)$ such that $\int_{0}^1 dr f(r)=1$ and $f$ is
a fixed point of $N$.

By (\ref{eqm1}), $m_1$ is a fixed point of $N$. Furthermore, if $V$ is
uniformly distributed over $\sun$ and independent of $\mathscr
{M}_{\infty}$, we know from
Corollary \ref{coruni} that~$\mathscr{M}_\infty(V)$ is the limit of
the uniformly integrable
martingale $\mathscr{M}^{(\alpha)}_t(V)$ (for any choice of $\alpha
\geq0$) and therefore
$\mathbb{E} [ \mathscr{M}_\infty(V) ]=1$. Hence,
\[
\int_{0}^1dr m_{1}(r)=\int_{0}^1dr \mathbb{E} [ \mathscr
{M}_{\infty }(e^{2i\pi r}) ]=\mathbb{E} [ \mathscr{M}_{\infty}(V) ]=1.
\]
We conclude that the function $f=m_{1}$ is the unique
nonnegative function in $\mathbb{L}^1(]0,1[,dr)$ such that $\int
_{0}^1 dr f(r)=1$ and $f$ is a fixed point of $N$.
On the other hand, elementary calculus shows that the function
$r\mapsto(r(1-r) )^{\beta^*}$
is also a fixed point of $N$. Indeed, noting that $2\beta^*/(\beta
^*+2)= 1-\beta^*$ and using
two integration by parts, we get
\[
\int_0^r \biggl(\frac{1-r}{1-u} \biggr)^{2+\beta^*} \biggl(\frac
{2}{1-u}-\frac{2\beta^*}{\beta^*+2} \biggr)\bigl (u(1-u)\bigr)^{\beta^*}\,
du=r^{1+\beta^*} (1-r)^{\beta^*}
\]
and similarly,
\[
\int_r^1 \biggl(\frac{r}{u} \biggr)^{2+\beta^*}\biggl (\frac
{2}{u}-\frac{2\beta^*}{\beta^*+2} \biggr)\bigl (u(1-u)\bigr)^{\beta^*}\,du=
r^{\beta^*} (1-r)^{\beta^*+1}.
\]
Therefore, the function
\[
e_{\beta^*}(r):= \frac{\Gamma(2+2\beta^*)}{\Gamma(1+\beta
^*)^2} \bigl(r(1-r) \bigr)^{\beta^*}
\]
is also a fixed point of $N$ such that $\int_{0}^1 dr e_{\beta
^*}(r)=1$. Consequently
we have $m_{1}(r) = e_{\beta^*}(r)$ a.e. The equality is in fact true
for every $r\in\,]0,1[$
since the integral equation (\ref{eqm1}) implies that $m_1$ is continuous
on $]0,1[$. This completes the
proof of Proposition \ref{moments} in the case $p=1$.

\subsection{\texorpdfstring{The case $p\geq2$}{The case p >= 2}}
From the H\"{o}lder inequality, and the case $p=1$, we have for every
integer $p\geq1$
and every $r\in\,]0,1[$,
%
%
\begin{equation}\label{holder}
m_{p}(r) \geq\biggl(\frac{\Gamma(2+2\beta^*)}{\Gamma(1+\beta
^*)^2} \biggr)^p\bigl(r(1-r)\bigr)^{p\beta^*}.
\end{equation}
We prove by induction on $k \geq1$, that for every $\varepsilon\in\,]0,1/2[$, there exists
a~constant $M_{\varepsilon,k}>0$ such that for
every $r\in\,]0,1[$,
%
%
\begin{equation}\label{rec}
m_{k}(r) \leq M_{\varepsilon,k}\bigl(r(1-r)\bigr)^{k\beta
^*-\varepsilon}.
\end{equation}
We assume that (\ref{rec}) holds for $k=1,2,\ldots,p-1$, and we prove
that this bound also
holds for $k=p$.

Similarly as in the case $p=1$, we can use the identity in distribution
(\ref{idloi}) to get
the following integral equation for the functions $m_{p}$:
%
%
\begin{eqnarray}
\label{eqnp}
m_{p}(r) &=& \int_{0}^r du \biggl(\frac{1-r}{1-u} \biggr)^{2+p\beta
^*} \biggl(\frac{2}{1-u}-\frac{2p\beta^*}{p\beta^*+2}
\biggr)m_{p}(u)\notag\\
&&{} + \int_{r}^1 du \biggl(\frac{r}{u} \biggr)^{2+p\beta^*}
\biggl(\frac{2}{u}-\frac{2p\beta^*}{p\beta^*+2}
\biggr)m_{p}(u)\nonumber
\\[-10pt]
\\[-10pt]
&&{} + 2\sum_{k=1}^{p-1} \pmatrix{p\cr k} \iint_{\mathscr{D}_{r}}
dr_{1}\,dr_{2}\nonumber\\
&&\hphantom{{}+ 2\sum_{k=1}^{p-1}}{}\times\frac{|r_{1}-r|^{1+(p-k)\beta^*}|r_{2}-r|^{1+k\beta
^*}}{|r_{1}-r_{2}|^{3+p\beta^*}}m_{k}(r_{1})m_{p-k}(r_{2}).
\nonumber
\end{eqnarray}
The derivation of (\ref{eqnp}) from (\ref{idloi}) is exactly similar
to that of (\ref{eqm1}), and we leave
details to the reader. Note that,
in contrast with the case $p=1$, we now get ``interaction terms''
involving the products $m_{k}(r)m_{p-k}(r)$. We start with some crude estimates.

\begin{lemma} \label{lemeap}For every $p \geq1$, the function $m_{p}$
is bounded over $]0,1[$. Moreover, for every $u,r\in\,]0,1/2[$, we have
%
%
\begin{equation} \label{utile}m_{p}(u) \leq
2^{p-1}\bigl(m_{p}(u+r)+m_{p}(r)\bigr).
\end{equation}
\end{lemma}

\begin{pf} For every $r,u\in\,]0,1[$, we set
\begin{eqnarray*}
g_{p,r}(u) &=&\mathbf{1}_{\{0< u < r\}} \biggl(\frac{1-r}{1-u}
\biggr)^{2+p\beta^*} \biggl(\frac{2}{1-u}-\frac{2p\beta^*}{p\beta
^*+2} \biggr)\\
&&{}
+\mathbf{1}_{\{r < u < 1\}} \biggl(\frac{r}{u} \biggr)^{2+\beta
^*}\biggl (\frac{2}{u}-\frac{2p\beta^*}{\beta^*+2} \biggr).
\end{eqnarray*}
From (\ref{eqnp}), we have
\[
m_p(r)\geq\int_{0}^1 du g_{p,r}(u)m_{p}(u).
\]
On the other hand, by using (\ref{idloi}) and the inequality $(a+b)^p
\leq2^{p-1}(a^p+b^p)$ for $a,b\geq0$,
we get
\begin{eqnarray*}
m_{p}(r)&\leq&2^{p-1} \bigl(\mathbb{E} [ M^{p\beta^*}\mathscr{M}'_{\infty
}(e^{2i\pi Y_1})^p ]
+ \mathbb{E} [ (1-M)^{p\beta^*}\mathscr{M}''_\infty(e^{2i\pi
Y_2})^p ] \bigr)\\
 &=&
2^{p-1} \int_{0}^1 du g_{p,r}(u)m_{p}(u),
\end{eqnarray*}
where the last equality again follows from calculations similar to
those leading to (\ref{eqm1}).

From the explicit form of the function $g_{p,r}$, we see that, for
every $\delta\in\,]0,1/2[$, there
exist positive constants $c_{\delta,p}$ and $C_{\delta,p}$ such that
for all $r \in\,]\delta,1-\delta[$,
%
%
\begin{equation} \label{bornitude}
c_{\delta,p} \int_{0}^1m_{p}(u)\,du
\leq m_{p}(r)\leq C_{\delta,p} \int_{0}^1m_{p}(u)\,du.
\end{equation}
If $U$ is uniformly distributed over $[0,1]$, Corollary \ref{coruni}
shows that\break
$\int_{0}^1m_{p}(u)\,du = \mathbb{E} [ \mathscr{M}_{\infty}(U)^p
]<\infty$. We
thus get that the function $m_{p}$ is bounded over every compact subset
of $]0,1[$.

To get information about the values of the
function $m_p$ in the neighborhood of $0$ (or of $1$), we use the
triangle inequality for figelas.\vadjust{\goodbreak}
Let $\alpha>0$. For every $r,u\in\,]0,1/2[$ and every $t\geq0$,
Proposition \ref{ineg} gives
\[
H_{S_\alpha(t)}(1,e^{2i\pi u})\leq H_{S_\alpha(t)}\bigl(1,e^{2i\pi(u+r)}\bigr)+
H_{S_\alpha(t)}\bigl(e^{2i\pi u}, e^{2i\pi(u+r)}\bigr).
\]
Furthermore, rotational invariance shows that the process
$(H_{S_\alpha(t)}(e^{2i\pi u},\break e^{2i\pi(u+r)}))_{t\geq0}$ has the
same distribution as the
process $(H_{S_\alpha(t)}(1,e^{2i\pi r}))_{t\geq0}$. We thus deduce from
Theorem \ref{tech}(iii) that
\[
\mathscr{M}_{\infty}(e^{2i\pi u})\leq\mathscr{M}_{\infty}\bigl(e^{2i\pi
(u+r)}\bigr) +\tilde{\mathscr M}_{\infty}(e^{2i\pi r}),
\]
where $\tilde{\mathscr M}_{\infty}(e^{2i\pi r})$ has the same
distribution as ${\mathscr M}_{\infty}(e^{2i\pi r})$.
The bound (\ref{utile}) now follows by using the inequality
$(a+b)^p\leq2^{p-1}(a^p+b^p)$ for
$a,b\geq0$.

Since we already know from (\ref{bornitude}) that the function $m_p$
is bounded over compact subsets of $]0,1[$,
and since $m_p(r)=m_p(1-r)$ by an obvious symmetry argument, the
bound (\ref{utile}) implies that $m_p$ is bounded over $]0,1[$.
\end{pf}

We come back to the proof of (\ref{rec}) with $k=p$. We fix
$\varepsilon\in\,]0,1/8[$. We start from the integral
equation (\ref{eqnp}) and first discuss the interaction terms. Fix $k
\in\{1,2, \ldots,p-1\}$ and set, for every
$r\in\,]0,1[$,
\[
T_{p,k}(r) = \iint_{\mathscr{D}_{r}} dr_{1}\,dr_{2}\frac
{|r_{1}-r|^{1+(p-k)\beta^*}|r_{2}-r|^{1+k\beta
^*}}{|r_{1}-r_{2}|^{3+p\beta^*}}m_{k}(r_{1})m_{p-k}(r_{2}).
\]
By the induction hypothesis, there exists a constant
$M_{p,k,\varepsilon}$ such that, for $r\in\,]0,1[$,
\[
T_{p,k}(r) \leq M_{p,k,\varepsilon}
\iint_{\mathscr{D}_{r}} dr_{1}\,dr_{2}\frac{|r_{1}-r|^{1+(p-k)\beta
^*}|r_{2}-r|^{1+k\beta^*}}{|r_{1}-r_{2}|^{3+p\beta^*}}r_{1}^{k\beta
^*-\varepsilon} r_{2}^{(p-k)\beta^*-\varepsilon}.
\]
Consider the integral over $[0,r]\times[r,1]$.
From the change of variables $r_1=rs_{1}$ and $r_2=rs_{2}$, we see that
this integral
is equal to
\begin{eqnarray*}
&&r^{p\beta^*+1-2\varepsilon}\int_{0}^1ds_{1}\int_{1}^{1/r}ds_{2}
\frac{(1-s_{1})^{1+(p-k)\beta^*}(s_{2}-1)^{1+k\beta
^*}}{(s_{2}-s_{1})^{3+p\beta^*}} s_1^{k\beta^*-\varepsilon
}s_{2}^{(p-k)\beta^*-\varepsilon}\\
&& \qquad \leq K r^{p\beta^* +1-2\varepsilon},
\end{eqnarray*}
where
\[
K = \int_{0}^1ds_{1}\int_{1}^\infty ds_{2}\frac
{(1-s_{1})^{1+(p-k)\beta^*}(s_{2}-1)^{1+k\beta
^*}}{(s_{2}-s_{1})^{3+p\beta^*}} s_{2}^{(p-k)\beta^*-\varepsilon} <
\infty.
\]
We get a similar bound for the integral over $[r,1]\times[0,r]$, and,
using the fact that $m_p(r)=m_p(1-r)$,
we conclude that the ``interaction terms'' in the integral equation
(\ref{eqnp}) are bounded above by a constant times $(r(1-r))^{p\beta
^*+1/2}$. By (\ref{holder}) these terms are negligible in comparison
with $m_{p}(r)$ when $r\to0$.

Thus for $r$ sufficiently close to $0$, say $0<r \leq r_{0} \leq1/4$,
we can write
\begin{eqnarray*}
m_{p}(r) &\leq&(1+\varepsilon)\int_{0}^r du \biggl(\frac
{1-r}{1-u} \biggr)^{2+p\beta^*} \biggl(\frac{2}{1-u}-\frac{2p\beta
^*}{p\beta^*+2} \biggr)m_{p}(u) \\
&&{} + (1+\varepsilon)\int_{r}^1 du \biggl(\frac{r}{u}
\biggr)^{2+p\beta^*} \biggl(\frac{2}{u}-\frac{2p\beta^*}{p\beta
^*+2} \biggr)m_{p}(u).
\end{eqnarray*}
The first term in the right-hand side is easily bounded by $3\int
_{0}^r du m_{p}(u)$, and we have, for $0<r\leq r_0$,
%
%
\begin{equation} \label{final}
m_{p}(r) \leq3\!\int_{0}^r\!du m_{p}(u) + (1+\varepsilon)\!\int_r^1\!du\biggl (\frac{r}{u}\biggr)^{2+p\beta^*}\!%
\biggl(\frac{2}{u}-\frac{2p\beta^*}{p\beta^*+2}\biggr)m_{p}(u).\hspace*{-30pt}
\end{equation}
However, by inequality (\ref{utile}), we have
for $0<r\leq r_0$,
%
%
\begin{equation}
\label{momenttechnical}
\int_{0}^{r} du m_{p}(u) \leq2^{p-1} \biggl(rm_{p}(r) + \int
_{r}^{2r}du m_{p}(u) \biggr).
\end{equation}
By (\ref{eqnp}), we have also
\[
m_p(r)\geq\int_{r}^{2r} du \biggl(\frac{r}{u} \biggr)^{2+p\beta
^*} \biggl(\frac{2}{u}-\frac{2p\beta^*}{p\beta^*+2} \biggr)m_{p}(u),
\]
and since $\frac{2}{u}$ tends to infinity as $u\to0$, this bound
shows that $\int_{r}^{2r}du m_{p}(u)$
is negligible in comparison with $m_p(r)$ when $r\to0$. Therefore,
from the bound
(\ref{momenttechnical}) and by choosing $r_0$ smaller if necessary, we
can assume that,
for $0<r\leq r_0$,
\[
3\int_{0}^r du m_{p}(u) \leq\biggl( 1 - \frac{1+\varepsilon
}{1+2\varepsilon} \biggr) m_{p}(r).
\]
By substituting this estimate in (\ref{final}), we get for $0<r\leq r_0$,
\[
m_{p}(r)\leq(1+2\varepsilon)\int_{r}^1 du \biggl(\frac{r}{u}
\biggr)^{2+p\beta^*} \biggl(\frac{2}{u}-\frac{2p\beta^*}{p\beta
^*+2} \biggr)m_{p}(u).
\]
Consequently, there exists a positive constant $K=K(r_{0},p,\varepsilon
)$ such that for $0<r\leq r_0$,
\[
\frac{m_{p}(r)}{r^{2+p\beta^*}} \leq K+
2(1+2\varepsilon)\int_r^{r_{0}} \frac{du}{u} \frac
{m_{p}(u)}{u^{2+p\beta^*}}.
\]
A straightforward application of Gronwall's lemma to the function $r\to
r^{-2-p\beta^*}\times m_p(r)$
gives for $0<r\leq r_0$,
\[
\frac{m_{p}(r)}{r^{2+p\beta^*}} \leq K \biggl(\frac{r_0}{r}
\biggr)^{2(1+2\varepsilon)},
\]
or equivalently
\[
m_p(r)\leq K r_0^{2(1+2\varepsilon)} r^{p\beta^*-4\varepsilon}.\vadjust{\goodbreak}
\]
Since $\varepsilon\in\,]0,1/8[$ was arbitrary, and since we have
$m_p(r)=m_p(1-r)$ for $r\in\,]0,1[$,
we have obtained the desired bound (\ref{rec}) at order $p$. This
completes the proof of
Proposition \ref{moments}.

\subsection{\texorpdfstring{Proof of Theorem \protect\ref{asympto}}{Proof of Theorem 1.1}}

The asymptotics in Theorem \ref{asympto} are consequences of the
more general results obtained in Corollary \ref{corofrag1}(ii)
and in Theorem \ref{tech}(iii), using also Remark \ref{Poissonconst}.
It remains to verify that the process $(\mathscr{M}_{\infty}(x), x
\in\sun)$ has a H\"{o}lder continuous modification. Let
$x$ and $y$ be two distinct points of $\sun\setminus\{0\}$, and let
$\alpha>0$.
By the triangle inequality in Proposition \ref{ineg}, we have for
every $t\geq0$,
\[
\bigl|H_{S_\alpha(t)}(1,x)-H_{S_\alpha(t)}(1,y)\bigr|\leq H_{S_\alpha(t)}(x,y),
\]
and $H_{S_\alpha(t)}(x,y)$ has the same distribution as
$H_{S_\alpha(t)}(1,x^{-1}y)$ by rotational invariance. We can let
$t\to\infty$ and
using Theorem \ref{tech}(iii), we get the
following stochastic inequality:
%
%
\begin{equation} \label{kolmo} |\mathscr{M}_{\infty}(e^{2i\pi
{r}})-\mathscr{M}_{\infty}(e^{2i\pi{s}})| \overset{(d)}{\leq}
\mathscr{M}_{\infty}\bigl(e^{2i\pi({r}-{s})}\bigr)
\end{equation}
for every $0\leq s< r<1$.

By Proposition \ref{moments}, we have then for every integer $p\geq1$
and every $0\leq s< r<1$,
\[
\mathbb{E} [ |\mathscr{M}_{\infty}(e^{2i\pi{r}})-\mathscr
{M}_{\infty }(e^{2i\pi{s}})|^p ] \leq M_{\varepsilon,p}
(r-s)^{p\beta
^*-\varepsilon}.
\]
Kolmogorov's continuity criterion (see \cite{RY99}, Theorem~I.2.1)
shows that the process
$(\mathscr{M}_{\infty}(x), x \in\sun)$ has a continuous modification,
which is even $(\beta^*-\varepsilon)$-H\"{o}lder continuous, for
every $\varepsilon>0$.

From now on, we only deal with the continuous modification
of the process $(\mathscr{M}_{\infty}(x), x \in\sun)$.
Recall the notation ${\mathcal T}_S$ for the plane tree associated with
a figela~$S$, and also recall that
$H_S$ corresponds to the graph distance on this tree.
One may ask about the convergence of the (suitably rescaled)
trees~${\mathcal T}_{S_{\alpha}(t)}$ in the sense of the Gromov--Hausdorff distance.
Recall the notation~$T_g$ for the $\mathbb R$-tree coded by a function
$g$ (see  Section~\ref{coding}).

\begin{conjecture*} Set ${g}_\infty(r)=\mathscr{M}_{\infty
}(e^{2i\pi r})$
for every $r\in[0,1]$. The convergence in distribution
\[
\bigl({\mathcal T}_{S_{\alpha}(t)}, t^{-\beta^*/\alpha}H_{S_{\alpha
}(t)} \bigr) \overset{(d)}{\underset{t \to\infty}{\longrightarrow
} } (T_{g_{\infty}},K_{\nu_{D}}(\alpha)\mathrm{d}_{g_{\infty}})
\]
holds in the sense of the Gromov--Hausdorff distance.
\end{conjecture*}

It would suffice to establish the following convergence in
distribution:
\[
\bigl( t^{-\beta^*/\alpha}H_{S_{\alpha}(t)}(1,e^{2i\pi r})
\bigr)_{r \in[0,1]} \overset{(d)}{\underset{t\to\infty}{\longrightarrow
}} (K_{\nu_{D}}(\alpha) \mathscr{M}_{\infty}(e^{2i\pi
r}) )_{r\in[0,1]}
\]
in the Skorokhod sense (the mapping $r\mapsto H_{S_{\alpha
}(t)}(1,e^{2i\pi r})$ is not defined when~$e^{2i\pi r}$
is a foot of $S_\alpha(t)$, but we can choose a suitable convention so
that\vadjust{\goodbreak} this mapping
is defined and c\`{a}dl\`{a}g over $[0,1]$).
Proving that this convergence holds would require more information
about the process $(H_{S_{\alpha}(t)}(1,x))_{x\in\sun, t\geq0}$.

\section{Identifying the limiting lamination}\label{sec5}

\subsection{Preliminaries}

The next proposition is the first step toward the proof
of Theorem \ref{codingL}. We recall the notation introduced at the
beginning of Section \ref{fragsepuni}: $a=e^{2i\pi U_1}$
and $b=e^{2i\pi U_2}$ are the feet of the first chord, with
$0<U_1<U_2<1$, and $M=1-(U_2-U_1)$.

\begin{proposition} \label{chord}
Conditionally on the pair $(U_1,U_2)$, we have
\[
\bigl(\mathscr{M}_\infty\bigl(e^{2i\pi(U_1+(U_2-U_1)r)}\bigr)-\mathscr
{M}_\infty(e^{2i\pi U_1}) \bigr)_{r\in[0,1]} \overset{(d)}{=}
\bigl( (1-M)^{\beta^*} \tilde{\mathscr M}_\infty(e^{2i\pi r})
\bigr)_{r\in[0,1]},
\]
where $\tilde{\mathscr M}_\infty$ is copy of $\mathscr{M}_\infty$
independent of $M$. Moreover, we have
\[
\mathscr{M}_\infty(e^{2i\pi U_1})>0 , \qquad  \mbox{a.s.}
\]
\end{proposition}

\begin{pf} This is essentially a consequence of Lemma \ref{keytool}. Fix
$\alpha>0$ and $r\in\,]0,1[$. Using the notation introduced before
this lemma, we have on the event $\{U_1<r<U_2\}$, for every $t\geq0$,
\[
H_{S_\alpha(\tau+t)}(1, e^{2i\pi r})= 1 + H_{S_\alpha^{(R')}(\tau
+t)}(1, e^{2i\pi U_1})+ H_{S_\alpha^{(R'')}(\tau+t)}(e^{2i\pi
U_1},e^{2i\pi r}).
\]
From Lemma \ref{keytool}, we now get on the event $\{U_1<r<U_2\}$ that
conditionally on $(U_1,U_2)$,
\begin{eqnarray*}
&&\bigl(H_{S_\alpha(\tau+t)}(1, e^{2i\pi r}) \bigr)_{t\geq0}\\
&& \qquad \overset{(d)}{=}
\bigl(1+ H_{S'_\alpha(M^\alpha t)}(1,\Psi_{a,b}(a)) + H_{S''_\alpha
((1-M)^\alpha t)}(1,\Phi_{a,b}(e^{2i\pi r}) )\bigr)_{t\geq0}.
\end{eqnarray*}
We multiply each side by $t^{-\beta^*/\alpha}$ and pass to the limit
$t\to\infty$, using Theo\-rem~\ref{tech}(iii), and we get with
an obvious notation that, on the event $\{U_1<r<U_2\}$ and
conditionally on $(U_1,U_2)$,
\[
\mathscr{M}_\infty(e^{2i\pi r}) \overset{(d)}{=} M^{\beta^*}
\mathscr{M}_\infty'(\Psi_{a,b}(a)) + (1-M)^{\beta^*} \mathscr
{M}_\infty''\bigl(e^{2i\pi\phi_{U_1,U_2}(r)}\bigr).
\]
This identity in distribution is immediately extended to a finite
number of values of $r$ by the same argument. Noting that
$\phi_{U_1,U_2}(U_1+(U_2-U_1)r)=r$, we thus get that, conditionally on
$(U_1,U_2)$,
\begin{eqnarray*}
&&\bigl(\mathscr{M}_\infty\bigl(e^{2i\pi(U_1+(U_2-U_1)r)}\bigr) \bigr)_{r\in
[0,1]}\\
&& \qquad  \overset{(d)}{=}
\bigl(M^{\beta^*}\mathscr{M}_\infty'(\Psi_{a,b}(a)) + (1-M)^{\beta
^*} \mathscr{M}_\infty''(e^{2i\pi r}) \bigr)_{r\in[0,1]}.
\end{eqnarray*}
In particular $\mathscr{M}_\infty(e^{2i\pi U_1})\overset{(d)}{=}
M^{\beta^*}\mathscr{M}_\infty'(\Psi_{a,b}(a))>0$ a.s. by Theorem
\ref{tech}(ii), and
the identity in distribution of the proposition also follows from the
previous\vadjust{\goodbreak} display.
\end{pf}

Recall the notation $S(\infty)$, $S^*(\infty)$ from the end of
Section~\ref{sec2}.

\begin{lemma}
\label{nochordfrom1}
For every $x\in\sun$, $\mathbb{P} [\exists y\in\sun\setminus\{x\}
\dvtx (x,y)\in S^*(\infty) ]=0$.
\end{lemma}

\begin{pf} Let $\varepsilon>0$. It is enough to prove that, for every
$x\in\sun$,
\[
\mathbb{P} [\exists y\in\sun\dvtx |y-x|>\varepsilon\mbox{ and
}(x,y)\in S^*(\infty) ]=0.
\]
Thanks to rotational invariance, this will follow if we can verify that
\[
\mathbb{E} \biggl[ \int\mathrm{m}(dx) {\mathbf1}_{\{ \exists y\in\sun\dvtx
|y-x|>\varepsilon\mbox{ and }(x,y)\in S^*(\infty)\}} \biggr]=0.
\]
Note that if $(x,y)\in S^*(\infty)$ the chord $[xy]$ does not cross
any of the (other) chords of
$S(\infty)$.

We can find an
integer $n$ (depending on $\varepsilon$) and $n$ points $z_1,\ldots
,z_n$ of $\sun$
such that the following holds. Whenever $x,y\in\sun$ are such that
$|y-x|>\varepsilon$, there exists
an index $j\in\{1,\ldots,n\}$ such that $z_j$ belongs to one of the
two open subarcs
with endpoints $x$ and $y$, and $-z_j$ belongs to the other subarc.
If we assume also that $(x,y)\in S^*(\infty)$, it follows that $x$
belongs to
the boundary of a fragment of $S_0(t)$ separating $z_j$ from $-z_j$,
for every
$t\geq0$.

Thanks to these observations, we have for every $t\geq0$,
\[
\int\mathrm{m}(dx) {\mathbf1}_{\{
\exists y\in\sun\dvtx |y-x|>\varepsilon\mathrm{\ and\ }(x,y)\in S^*(\infty)\}}
\leq\sum_{j=1}^n \biggl(\sum_i \mathrm{m}(R^{z_j,-z_j}_i(S_0(t))
)\biggr)
\]
with the notation introduced in Definition \ref{sdc}. From Theorem
\ref{tech}(iv) and the fact
that $\kappa_{\nu_D}(1)>0$,
the right-hand side tends to $0$ almost surely as $t\to\infty$, which
completes the proof.
\end{pf}

Recall our notation
$g_\infty(r)=\mathscr{M}_\infty(e^{2i\pi r})$ for every $r\in
[0,1]$. Notice that
$g=g_\infty$ satisfies the assumptions of  Section~\ref{coding}.

\begin{cor}
\label{codinglami}
Almost surely, for every $r,s\in[0,1]$ such that $\{e^{2i\pi
r},\allowbreak e^{2i\pi s}\}\in S(\infty)$, we have $r\overset{g_\infty}{\approx
} s$.
\end{cor}

\begin{pf} If $c,d$ are two distinct points of $\sun\setminus\{1\}$,
write $\op{Arc}^*(c,d)$ for the
open subarc of $\sun$ with endpoints $c$ and $d$
not containing $1$. As an immediate consequence of Proposition \ref
{chord}, we have $\mathscr{M}_\infty(x)\geq\mathscr{M}_\infty(a)
=\mathscr{M}_\infty(b)>0$, for every $x\in\op{Arc}^*(a,b)$. This
property is easily extended by induction (using Lemma \ref{keytool}
once again)
to any chord appearing in the figela process. We have almost surely for
every $\{c,d\}\in S(\infty)$,
%
%
\begin{equation}
\label{positiv}
\mathscr{M}_\infty(x)\geq\mathscr{M}_\infty(c)
=\mathscr{M}_\infty(d)>0 \qquad  \mbox{for every }x\in\op{Arc}^*(c,d).
\end{equation}

We can in fact replace the weak inequality $\mathscr{M}_\infty(x)\geq
\mathscr{M}_\infty(c)$ by a strict one. To see this, we first note that, by
Lemma \ref{nochordfrom1}, $1$ is not an endpoint of a
(nondegenerate) chord of $S^*(\infty)$.\vadjust{\goodbreak} By an easy argument, this
implies that almost surely, for every
$\varepsilon>0$, there exist $r\in\,]-\varepsilon,0[$ and $s\in\,]0,\varepsilon[ $ such
that the chord $[e^{2i\pi r}e^{2i\pi s}]$
belongs to $S(\infty)$. It follows that
\[
\bigcup_{\{c,d\}\in S(\infty)} \op{Arc}^*(c,d)= \sun\setminus\{1\}
 \qquad  \mbox{a.s.}
\]
From (\ref{positiv}) we now get that $\mathscr{M}_\infty(x)>0$, for
every $x\in\sun\setminus\{1\}$, a.s.

We can apply this property to the process $\tilde{\mathscr M}_\infty$
in Proposition \ref{chord}, and we get
that $\mathscr{M}_\infty(x)>\mathscr{M}_\infty(a)=\mathscr
{M}_\infty(b)$, for every $x\in\op{Arc}^*(a,b)$, a.s. Again, this property
of the first chord is easily extended by induction to any chord in the
figela process, and we obtain that,
almost surely for every $\{c,d\}\in S(\infty)$,
%
%
\begin{equation}
\label{positivstrict}
\mathscr{M}_\infty(x)> \mathscr{M}_\infty(c)
=\mathscr{M}_\infty(d) \qquad  \mbox{for every }x\in\op{Arc}^*(c,d).
\end{equation}
The statement of Corollary \ref{codinglami} now follows from the
definition of $\overset{g_\infty}{\approx}$.
\end{pf}

If $(x,y)\in S^*(\infty)$ we can write $(x,y)=\lim(x_n,y_n)$ where\vspace*{1pt} $\{
x_n,y_n\}\in S(\infty)$ for every $n$.
Write $x=e^{2i\pi r}$, $y=e^{2i\pi s}$ and $x_n=e^{2i\pi r_n}$,
$y_n=e^{2i\pi s_n}$,
where~$r,s,\allowbreak r_n,s_n\in[0,1]$. By Corollary \ref{codinglami}, we have
$r_n\overset{g_\infty}{\approx} s_n$
for every $n$.
Since the graph of the relation $\overset{g_\infty}{\approx} $ is
closed, it follows that $r\overset{g_\infty}{\approx} s$.
We have thus proved that
\[
L_\infty\subset L_{g_\infty}.
\]
The reverse inclusion will be proved in the next subsection.

\subsection{Maximality of the limiting lamination}

The proof of Theorem \ref{codingL} will be completed thanks to the
following proposition.

\begin{proposition}\label{max} Almost surely, $L_{\infty}$ is a
maximal lamination of $\overline{\mathbb{D}}.$
\end{proposition}

\begin{figure}

\includegraphics{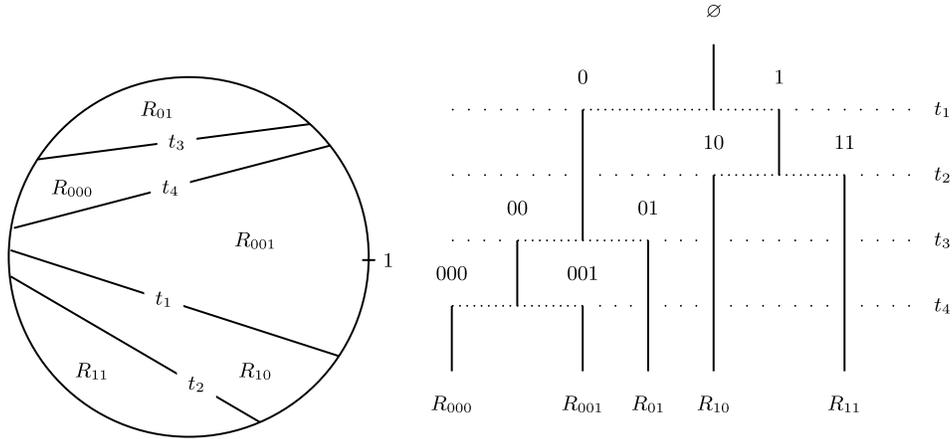}

\caption{Chords are represented on the left-hand side with their
respective creation times. On the right-hand side the genealogical tree
of the fragments $R_{000}, \ldots , R_{11}$ present at time~$t_{4}$.}
\label{fig5}
\end{figure}

Before proving Proposition~\ref{max}, we need to
establish a preliminary lemma. This lemma is concerned
with the genealogical tree of fragments appearing in the figela
process, which we construct as follows. We consider
the
fragments created by $S_{0}(t)$ as time increases. The first fragment
is $R_\varnothing=\overline{\mathbb{D}}$. At the exponential time
$\tau$, the first chord splits $\overline{\mathbb{D}}$ into two
fragments, which are viewed as the offspring of $\varnothing$.
We then order these fragments in a random way: with probability $1/2$,
we call $R_0$ the
fragment with the largest mass and $R_1$ the other one, and with
probability $1/2$ we do the contrary.
We then iterate this device. Then each fragment that appears in the
figela process is labeled by an element of
the infinite binary tree
\[
\mathbb{T}=\bigcup_{n \geq0} \{0,1\}^n.
\]
For every integer $n$, we also set
\[
\mathbb{T}_n:=\bigcup_{k=0}^n \{0,1\}^k.\vadjust{\goodbreak}
\]
At every time $t$, we have a (finite) binary tree corresponding to the
genealogy of the fragments
present at time $t$. See Figure~\ref{fig5}.

%

If $R$ is a fragment, we call \textit{end} of $R$ any connected
component of
$R\cap\sun$. We denote the number of ends of a fragment $R$ by $\mathrm{e}(R).$
For reasons that will be explained later, the full disk $\overline
{\mathbb{D}}$ is viewed as a
fragment with $0$ end.

\begin{lemma} \label{3bouts} In the (infinite) genealogical tree of
fragments, almost surely, there is no ray along which all fragments
have eventually strictly more than~3 ends.
\end{lemma}

\begin{pf}
Let $n\geq0$ and $u\in\{0,1\}^n$, and consider the fragment $R_{u}$.
Let~$y_u$ be
one endpoint (chosen at random) of the first chord that will fall
inside~$R_u$.
Note that, conditionally on $R_u$, $y_u$ is uniformly distributed over
$R_u\cap\sun$. Let $\varphi_u\dvtx[0,\mathrm{m}(R_u)[\longrightarrow
\overline{R_u\cap\sun}$ be defined by requiring that the measure of the
intersection of $R_u\cap\sun$ with the arc
(in counterclockwise order) between~$y_u$ and~$\varphi_u(t)$ is equal
to $t$, for every $t\in[0,\mathrm{m}(R_u)[$.
This definition is unambiguous if we also impose that $\varphi_u$ is
right-continuous. Then
$\varphi_u$ has exactly~$\mathrm{e}(R_u)$ discontinuity times
corresponding to the chords that lie in the boundary of~$R_u$ (indeed the left and right limits of $\varphi_u$ at a
discontinuity time are
the endpoints of a chord adjacent to $R_u$).
We claim that, conditionally given $(\mathrm{m}(R_u),\mathrm{e}(R_u))$, the
set of discontinuity times of $\varphi_u$
is distributed as the collection of $\mathrm{e}(R_u)$ independent points
chosen uniformly over $[0,\mathrm{m}(R_u)[$.

This claim can be checked by induction on $n$. For $n=0$ there is
nothing to prove. Assume that
the claim holds up to order $n$. Recalling that~$y_u$ is one endpoint
of the first chord that will fall in $R_u$, the other endpoint~$z_u$
will be chosen uniformly\vadjust{\goodbreak} over $R_u\cap\sun$, so that $\varphi_u^{-1}(z_u)$
will be uniform over $[0,\mathrm{m}(R_u)[$. We have then $\mathrm
{e}(R_{u0})=K+1$ and $\mathrm{e}(R_{u1})=\mathrm{e}(R_u)+1-K$
(or the contrary with probability $1/2$),
where $K$ is the number of discontinuity times of $\varphi_u$ in
$[0,\varphi_u^{-1}(z_u)]$.
Using our induction hypothesis, we see that
conditionally on $K$ and on $\varphi_u^{-1}(z_u)$ the latter
discontinuity times are independent and uniformly distributed over
$[0,\varphi_u^{-1}(z_u)]$, and that a similar property holds for the
discontinuity times that belong
to $[\varphi_u^{-1}(Z_u),\mathrm{m}(R_u)]$. It follows that the desired
property will still hold at order $n+1$.

The preceding arguments also show that, conditionally on $R_u$, $\mathrm
{e}(R_{u0})$ is distributed as
$K+1$, where $K$ is obtained by throwing $\mathrm{e}(R_u)+1$ uniform
random variables in $[0,\mathrm{m}(R_u)]$
and counting how many among the $\mathrm{e}(R_u)$ first ones are smaller
than the last one.
By an obvious symmetry argument, we have, for any integers $p\geq0$,
$k\in\{0,\ldots,p\}$ and any $a\in\,]0,1]$,
\[
\mathbb{P} [\mathrm{e}(R_{u0})=k+1 | \mathrm{e}(R_{u})=p, \mathrm{m}(R_u)=a]
= \frac{1}{p+1}.
\]
Notice that the preceding conditional probability does not depend on
$a$, which could have been seen
from a scaling argument.

Modulo some technical details that are left to the reader, we get that
the distribution of the tree-indexed process
$(\mathrm{e}(R_u),u\in\mathbb{T})$ can be described as follows.
We start with $\mathrm{e}(\varnothing)=0$, and we then proceed by
induction on~$n$
to define~$\mathrm{e}(R_u)$ for every $u\in\{0,1\}^n$. To this end, given
the values of~$\mathrm{e}(R_u)$ for $u\in\mathbb{T}_n$,
we choose independently for every $v\in\{0,1\}^{n}$ a random variable~$k_v$ uniform over
$\{0,\ldots,\mathrm{e}(R_v)\}$, and we set $\mathrm{e}(R_{v0})=k_v+1$, $\mathrm
{e}(R_{v1})=\mathrm{e}(R_v)-k_v+1$.

Consider a tree-indexed process $(\mathrm{f}_u,u\in\mathbb{T})$ that
evolves according to the preceding
rules but starts with $\mathrm{f}_\varnothing=4$ [instead of $\mathrm
{e}(R_\varnothing)=0$]. In order to get the statement
of the lemma, it is enough to prove that almost surely, there is no
infinite ray starting from the root along which all the values of~$\mathrm{f}_u$ are strictly larger than~$3$.
Consider a fixed infinite ray in the tree, say
$\varnothing,0,00,\allowbreak 000,\ldots$ and let $X_{0}=\mathrm{f}_{\varnothing}$,
$X_1=\mathrm{f}_0$,
$X_2=\mathrm{f}_{00},\ldots$ be the
values of our process along the ray. Note that $(X_{n})_{n\geq0}$ is a
Markov chain
with values in~$\mathbb{N}$, with transition kernel given by
\[
q_{k\ell}= \frac{1}{k+1} \mathbf{1}_{\{1\leq l\leq k+1\}}
\]
for every $k,\ell\geq1$. Write $(\mathcal{F}_{n})_{n\geq0}$ for the
filtration generated by the process $(X_{n})_{n\geq0}$. We have
\[
\mathbb{E} [ X_{n+1}|\mathcal{F}_{n} ] = \frac{1}{X_{n}+1} \bigl(1+ 2
+ \cdots+(X_{n}+1)\bigr)
= \frac{X_{n}}{2}+1.
\]
Hence $M_{n}= 2^n(X_{n}-2)$ is a martingale starting from $2$. For $i
\geq1$ we let $T_{i}$ be the stopping time $T_{i}=\inf\{n \geq0\dvtx
X_{n}=i\},$ and $T= T_{1}\wedge T_{2}\wedge T_{3}$.

Note that $\mathbb{P} [X_{k} \geq4, \mbox{ for every } 0\leq k \leq
n ] = \mathbb{P} [T> n ]$,
and that the preceding discussion applies to the values of $f_u$ along any
infinite ray starting from the\vadjust{\goodbreak} root.

By the stopping theorem applied to the martingale $(M_{n})_{n\geq0}$,
we obtain for every $n \geq0$,
\begin{eqnarray*}
2 &=&\mathbb{E} [ M_{n \wedge T} ]\\ &=&\mathbb{E} \bigl[ -2^{T_{1}}\mathbf
{1}_{\{T_{1}=T \leq n\}} \bigr]
+ 0 +
\mathbb{E} \bigl[ 2^{T_{3}}\mathbf{1}_{\{T_{3}=T \leq n\}} \bigr] +\mathbb{E}
\bigl[ 2^n(X_{n}-2)\mathbf {1}_{\{T>n\}} \bigr].
\end{eqnarray*}
From the transition kernel of the Markov chain $(X_{n})_{n\geq0}$ it
is easy to check that for every $k \geq1$,
$ \mathbb{P} [T_{1}=T=k ]= \mathbb{P} [T_{2}=T=k ]=\mathbb{P}
[T_{3}=T=k ].$ Hence, the equality
in the last display becomes
\[
2= \mathbb{E} \bigl[ 2^n(X_{n}-2)\mathbf{1}_{\{T>n\}} \bigr],
\]
or equivalently
\[
2= 2^n \mathbb{P} [T>n ] \mathbb{E} [ X_{n}-2\mid T>n ].
\]
Since obviously $\mathbb{E} [ X_{n}-2\mid T>n ]\geq2$, we get $2^n
\mathbb{P} [T>n ]
\leq1$.

For every $u=(u_1,\ldots,u_n)\in\{0,1\}^n$, and every $j\in\{
0,1,\ldots,n\}$, set $[u]_j=(u_1,\ldots,u_j)$, and
if $j\geq1$, also set
$[u]_j^*=(u_1,\ldots,u_{j-1},1-u_{j})$. Let
\[
G_n=\bigl\{u\in\{0,1\}^n\dvtx \mathrm{f}_{[u]_j}\geq4 , \forall j\in\{
0,1,\ldots,n\}\bigr\}.
\]
Clearly
%
%
\begin{equation}
\label{bout1}
\mathbb{E} [ \# G_n ]=2^n \mathbb{P} [T>n ] \leq1.
\end{equation}
In order to get the statement of the lemma, it is enough to verify that
$\mathbb{P} [\# G_n\geq1 ]\longrightarrow0$
as $n\to\infty$. Note that the sequence $\mathbb{P} [\# G_n\geq1 ]$ is
monotone nonincreasing.
We argue by
contradiction and assume that there exists $\eta>0$ such that $\mathbb
{P} [\# G_n\geq1 ]\geq\eta$
for every $n\geq1$. By a simple
coupling argument, the same lower bound will remain valid if we start
the tree-indexed process
with $\mathrm{f}_\varnothing= m$, for any $m\geq4$, instead of $\mathrm
{f}_\varnothing= 4$.

Fix $\varepsilon\in\,]0,\eta[$, and choose an integer $\ell\geq1$
such that $1/\ell\leq\varepsilon/2$. Choose another integer $k\geq
1$ such that,
if $B_{k,\eta}$ denotes a binomial ${\mathcal B}(k,\eta)$ random
variable, we have $\mathbb{P} [B_{k,\eta}\leq\ell ] \leq
\varepsilon/2$.
Finally set
\begin{eqnarray*}
G'_n&=&\bigl\{u\in G_n\dvtx \#\bigl\{j\in\{0,1,\ldots,n-1\}\dvtx\mathrm
{f}_{[u]_{j+1}}\leq
\mathrm{f}_{[u]_j}-2\bigr\}\leq k\bigr\},\\
G''_n &=& G_n\setminus G'_n.
\end{eqnarray*}

We first evaluate $\mathbb{P} [G''_n\not= \varnothing ]$. We have
\[
\mathbb{P} [G''_n\not= \varnothing ] \leq\mathbb{P} [G''_n\not=
\varnothing, \# G_n\leq\ell ] + \mathbb{P} [\# G_n>\ell ].
\]
By (\ref{bout1}) and our choice of $\ell$, we have $\mathbb{P} [\#
G_n>\ell ]\leq\ell^{-1}\mathbb{E} [ \#G_n ]\leq\varepsilon/2$. On
the other hand,
\begin{eqnarray*}
\mathbb{P} [G''_n\not= \varnothing, \# G_n\leq\ell ] &\leq
&\mathbb{E} \bigl[ \# G''_n {\mathbf1}_{\{ \# G_n\leq\ell\}} \bigr]\\
&=&\sum_{u\in\{0,1\}^n} \mathbb{P} [u\in G''_n, \# G_n\leq\ell ]\\
&=&\sum_{u\in\{0,1\}^n} \mathbb{P} [u\in G''_n ] \mathbb{P} [\#
G_n\leq\ell| u\in G''_n ].
\end{eqnarray*}
Fix $u\in\{0,1\}^n$. We argue conditionally on the values
of $\mathrm{f}_{[u]_{j}}$ for $0\leq j\leq n$, and note that the values of
$\mathrm{f}_{[u]_{j+1}^*}$
for $0\leq j\leq n-1$
are then also determined by the condition $\mathrm{f}_{[u]_{j+1}}+\mathrm
{f}_{[u]_{j+1}^*} = \mathrm{f}_{[u]_j} +2$.
Moreover, on the event $\{u\in G''_n\}$, there are at least $k$ values
of $j\in\{0,1,\ldots,n-1\}$
such that \mbox{$\mathrm{f}_{[u]_{j+1}}\leq\mathrm{f}_{[u]_j}-2$}.
For these values of $j$, we must have
$\mathrm{f}_{[u]_{j+1}^*}\geq4$. Furthermore, for each such value of $j$,
there is (conditional) probability at least $\eta$ that one of the
descendants of $[u]^*_{j+1}$ at generation $n$, say $v$, is
such that $\mathrm{f}_{[v]_i}\geq4$ for every $i\in\{j+1,\ldots,n\}$,
and consequently $v\in G_n$.
Summarizing, we see that conditionally on the event $\{u\in G''_n\}$,
$\# G_n$ is bounded
below in distribution by a binomial ${\mathcal B}(k,\eta)$ random
variable. Hence, using
our choice of $k$,
\[
\mathbb{P} [\# G_n\leq\ell| u\in G''_n ]\leq\mathbb{P}
[B_{k,\eta}\leq\ell ]
\leq\frac{\varepsilon}{2}.\vspace*{-1pt}
\]
We thus get
\[
\mathbb{P} [G''_n\not= \varnothing, \# G_n\leq\ell ]
\leq\frac{\varepsilon}{2} \sum_{u\in\{0,1\}^n} \mathbb{P} [u\in
G''_n ] =
\frac{\varepsilon}{2}\mathbb{E} [ \# G''_n ] \leq\frac{\varepsilon}{2},\vspace*{-1pt}
\]
by (\ref{bout1}). It follows that
\[
\limsup_{n\to\infty} \mathbb{P} [G''_n\not=\varnothing ] \leq
\varepsilon.\vspace*{-1pt}
\]

We will now verify that $\mathbb{P} [G'_n\not=\varnothing ]$ tends
to $0$
as $n\to\infty$. Since $\varepsilon<\eta$, this will give a
contradiction with our
assumption $\mathbb{P} [\# G_n\geq1 ]\geq\eta$
for every $n\geq1$, thus completing the proof.
We in fact show that $\mathbb{E} [ \#G'_n ]$ tends to $0$
as $n\to\infty$. To this end,
we first write, for $n\geq k$,
%
%
\begin{eqnarray}
\label{bout2}
\mathbb{E} [ \#G'_n ]
&=&2^n \mathbb{P} [T>n, \#\bigl\{j\in\{0,\ldots
,n-1\}\dvtx X_{j+1}\leq X_j -2\bigr\}\leq k ]\nonumber
\\[-1pt]
&\leq&2^n n^k \sup_{A\subset\{0,1,\ldots,n-1\},\# A= k} \mathbb{P}
[X_{j+1}\geq(X_j-1)\vee4 ,\\[-1pt]
&&\hphantom{\;2^n n^k \sup_{A\subset\{0,1,\ldots,n-1\},\# A= k} \mathbb{P}} \forall j\in\{0,1,\ldots,n-1\} \setminus
A ].
\nonumber\vspace*{-1pt}
\end{eqnarray}
We thus need to bound the quantity
\[
\mathbb{P} [X_{j+1}\geq(X_j-1)\vee4 , \forall j\in\{0,1,\ldots
,n-1\} \setminus A ],\vspace*{-1pt}
\]
for every choice of $A\subset\{0,1,\ldots,n-1\}$ such that $\# A= k$.
For every subset~$A$ of $\{0,1,\ldots,n-1\}$, we set
\[
N^A_n=\#\bigl\{j\in\{0,1,\ldots,n-1\}\setminus A\dvtx X_j =5\bigr\}.\vspace*{-1pt}
\]
With a slight abuse of notation, write $\mathbb{P}_i$ for a
probability measure under
which the Markov chain $X$ starts from $i$. We prove by induction on $n$
that for every choice of $A\subset\{0,1,\ldots,n-1\}$ and
$m\in\{0,1,\ldots,n-\# A\}$, we have
for every $i\geq1$,
%
%
\begin{eqnarray}
\label{bout3}
&&{\mathbb P}_i [X_{j+1}\geq(X_j-1)\vee4 , \forall j\in\{
0,1,\ldots,n-1\}\setminus A ; N^A_n=m ]\nonumber
\\[-10pt]
\\[-10pt]
&& \qquad \leq\biggl(\frac{1}{2} \biggr)^m
\biggl(\frac{3}{7} \biggr)^{n-m-\#A}.
\nonumber\vadjust{\goodbreak}
\end{eqnarray}
If $n=0$ (then necessarily $m=0$ and $A=\varnothing$) there is nothing
to prove. Assume that the desired bound holds at order $n-1$.
In order to prove that it holds at order $n$, we apply the Markov
property at time $1$. We need to
distinguish three cases.

If $0\in A$, then the left-hand side of (\ref{bout3}) is equal to
\[
\sum_{i'} q_{ii'} \mathbb{P}_{i'} [X_{j+1}\geq(X_j-1)\vee4
, \forall j\in\{0,1,\ldots,n-2\}\setminus A' ;
N^{A'}_{n-1}=m ],
\]
where $A'=\{j-1\dvtx j\in A,j>0\}$. Since $\# A'=\# A -1$ in that case, an
application of the induction hypothesis gives the result.

If $0\notin A$ and $i\not=5$, then the left-hand side of (\ref
{bout3}) is equal to
\begin{eqnarray*}
&&\sum_{i'\geq(i-1)\vee4} q_{ii'} \mathbb{P}_{i'} [X_{j+1}\geq
(X_j-1)\vee4 , \forall j\in\{0,1,\ldots,n-2\}\setminus A' ;
N^{A'}_{n-1}=m ]\\
&& \qquad \leq\sum_{i'\geq(i-1)\vee4} q_{ii'}\biggl (\frac{1}{2} \biggr)^m
\biggl(\frac{3}{7} \biggr)^{n-1-m-\#A'},
\end{eqnarray*}
and we just have to observe that
\[
\sum_{i'\geq(i-1)\vee4} q_{ii'} \leq\frac{3}{7}
\]
when $i\not=5$.

Finally, if $0\notin A$ and $i=5$, the left-hand side of (\ref{bout3})
is equal to
\begin{eqnarray*}
&&\sum_{i'\geq4} q_{5i'} \mathbb{P}_{i'} [X_{j+1}\geq
(X_j-1)\vee4 , \forall j\in\{0,1,\ldots,n-2\}\setminus A' ;
N^{A'}_{n-1}=m-1 ]\\
&& \qquad \leq\biggl(\sum_{i'\geq4} q_{5i'} \biggr)\biggl (\frac
{1}{2} \biggr)^{m-1}
\biggl(\frac{3}{7} \biggr)^{(n-1)-(m-1)-\#A'}\\
&& \qquad = \biggl(\frac{1}{2} \biggr)^m
\biggl(\frac{3}{7} \biggr)^{n-m-\#A},
\end{eqnarray*}
using the fact that $\sum_{i'\geq4} q_{5i'} =1/2$. This completes the
proof of (\ref{bout3}).

Fix $\delta\in\,]0,1[$. By summing over possible values of $m$, we
get for $n$ large, for all choices of
$A\subset\{0,1,\ldots,n-1\}$ such that $\# A= k$
\begin{eqnarray*}
&&\mathbb{P} [X_{j+1}\geq(X_j-1)\vee4 , \forall j\in\{0,1,\ldots
,n-1\} \setminus A ]\\
&& \qquad \leq\sum_{m=0}^{n-\lfloor\delta n\rfloor} \biggl(\frac
{1}{2} \biggr)^m
\biggl(\frac{3}{7} \biggr)^{n-m-k} + \mathbb{P} [N^A_n > n-\lfloor\delta
n\rfloor ]\\
&& \qquad \leq n \biggl(\frac{1}{2} \biggr)^{n-\lfloor\delta n\rfloor
} \biggl(\frac{3}{7} \biggr)^{\lfloor\delta n\rfloor-k}
+ \mathbb{P} [N^A_n > n-\lfloor\delta n\rfloor ].
\end{eqnarray*}
Note that $N^A_n\leq N^\varnothing_n$.
Crude estimates, using the fact that $\sup_{i\geq1} q_{i5} =1/5$,
show that we can fix $\delta$ such that
\[
2^n n^{k+1} \mathbb{P} [N^\varnothing_n > n -\lfloor\delta n\rfloor ]
\underset{n\to\infty}{\longrightarrow} 0.
\]
It then follows that the right-hand side of (\ref{bout2}) tends to $0$
as $n\to\infty$, which completes the proof.
\end{pf}

\begin{rek} For every integer $k\geq1$, let $(\mathrm
{f}^{(k)}_u,u\in\mathbb{T})$ be a tree-indexed process that evolves
according to the same rules
as $(\mathrm{e}(R_u),u\in\mathbb{T})$ but starts with $\mathrm
{f}^{(k)}_\varnothing=k$. Let $p_k$ be the probability that there
exists no infinite ray starting from $\varnothing$ along which
all labels $\mathrm{f}^{(k)}_u$ are strictly greater than $3$. By
conditioning on
the values of $\mathrm{f}^{(k)}_0$ and $\mathrm{f}^{(k)}_1$, we see that
$(p_k)_{k\geq1}$ satisfies
the properties
%
%
\begin{equation}
\label{recursiv}
 \quad \cases{
\displaystyle
p_{1}=p_{2}=p_{3}=1, \vspace*{2pt}\cr\displaystyle
p_{k} = \frac
{1}{k+1}(p_{1}p_{k+1}+ p_{2}p_{k}+ \cdots +p_{k}p_{2}+p_{k+1}p_{1}),&
\quad
if $k \geq4$.
}
\end{equation}
It follows that the values of $p_k$ for $k \geq5$ are determined recursively
from the value of $p_4$. Numerical simulations suggest that there
exists no sequence
$(p_k)_{k\geq1}$ satisfying (\ref{recursiv}) such that $p_4<1$ and
$0\leq p_k\leq1$
for every $k\geq1$. A~rigorous verification of this fact would provide an
alternative more analytic proof of Lemma \ref{3bouts}.
\end{rek}

\begin{pf*}{Proof of Proposition \ref{max}} First note that it
is easy to verify that $L_\infty\cap\sun$ is dense in
$\sun$, and thus $\sun\subset L_\infty$ since $L_\infty$ is closed.
We argue by contradiction and suppose that $L_\infty$ is not
a maximal lamination. Then there exists a (nondegenerate) chord $[xy]$
which is not contained in $L_\infty$ and is
such that $L_\infty\cup[xy]$ is still a lamination,
which implies that $]xy[$ does not intersect any chord of $S(\infty)$.
There is a unique infinite ray $\varnothing, \epsilon_1,
\epsilon_1\epsilon_2,\ldots$ in $\mathbb T$ such that $]xy[\subset
R_{\epsilon_1\cdots\epsilon_n}$ for every integer $n\geq0$.
We claim that, for all sufficiently large $n$, $R_{\epsilon_1\cdots
\epsilon_n}$ has at least $4$ ends. To see this, denote by $I^x_n$ the
end of $R_{\epsilon_1\cdots\epsilon_n}$ whose closure $\overline
{I^x_n}$ contains $x$, and define $I^y_n$ similarly.
Note that the maximal length of an end of a fragment at the $n$th generation
tends to $0$ a.s., and that this applies in particular to $I^x_n$ and
$I^y_n$. It follows that, almost surely for all $n$
sufficiently large, there is no chord of $S(\infty)$ between a point
of $\overline{I^x_n}$ and a point of $\overline{I^y_n}$ (otherwise, the
pair $(x,y)$ would be in $S^*(\infty)$ and the chord $[xy]$ would be
contained in $L_\infty$). Hence, for all
sufficiently large $n$, the boundary of $R_{\epsilon_1\cdots\epsilon
_n}$ contains at least $4$ different chords, and therefore
at least $4$ ends. This contradicts
Lemma \ref{3bouts}, and this contradiction completes the proof.
\end{pf*}

\begin{pf*}{Proof of Theorem \ref{codingL}} Since $L_\infty$ is
a maximal lamination and $L_\infty\subset L_{g_\infty}$,
we must have $L_\infty= L_{g_\infty}$ and in particular $L_{g_\infty
}$ is a maximal lamination. Thus, the
function $g_\infty$ must satisfy\vadjust{\goodbreak} the necessary and sufficient
condition for maximality given in
Proposition \ref{codinglamination}. Under this condition, however, the
relations $\overset{g_\infty}{\approx}$
and $\overset{g_\infty}{\sim}$ coincide. Recalling that ${\mathscr
M}_\infty(x)>0={\mathscr M}_\infty(1)$ for every
$x\in\sun\setminus\{1\}$, we see that property (\ref{maxipro})
written with $x=e^{2i\pi r}$ and
$y=e^{2i\pi s}$ is equivalent to saying that $r\overset{g_\infty
}{\sim} s$. Theorem \ref{codingL} then follows from the fact
that $L_\infty= L_{g_\infty}$.
\end{pf*}

\begin{rek}
\label{chordsinL}
It is not hard to see that $L_\infty$ has zero Lebesgue measure
a.s. (this follows from
the upper bound on the Hausdorff dimension proved in the next section).
By a simple argument,
it follows that a chord $[xy]$ is contained in $L_\infty$ if and only
if $x\overset{g_\infty}{\approx}y$,
and this condition is also equivalent to $(x,y)\in S^*(\infty)$.
\end{rek}

\section{\texorpdfstring{The Hausdorff dimension of $L_\infty$}{The Hausdorff dimension of L infinity}}\label{sec6}

In this section, we prove Theorem~\ref{Hausdim}.
We let ${\mathcal I}$ be the countable set of all pairs $(I,J)$ where
$I$ and $J$
are two disjoint closed subarcs of $\sun$ with nonempty interior and
endpoints of the form
$\exp(2i\pi r)$ with rational $r$. For each $(I,J)\in{\mathcal I}$,
we set
\[
L_{(I,J)}=\bigcup_{(y,z)\in S^*(\infty)\cap(I\times J)} [yz] \subset
L_\infty.
\]
Clearly,
%
%
\begin{equation}
\label{HausdorffTech}
\op{dim} L_\infty=\sup_{(I,J)\in{\mathcal I}} \op{dim}\bigl(L_{(I,J)}\bigr).
\end{equation}

\subsection*{Upper bound} We prove that, for every $(I,J)\in{\mathcal I}$,
\[
\op{dim}\bigl(L_{(I,J)}\bigr) \leq\frac{\sqrt{17}-1}{2} = \beta^* +1
\qquad\mbox{a.s.}
\]
By rotational invariance, we may assume without loss of generality that
$1\notin I\cup J$.
We pick a point $x\in\sun\setminus(I\cup J)$ such that $1$ and $x$
belong to
different components of $ \sun\setminus(I\cup J)$. We also fix
$\gamma> \beta^* +1$ and
set $\beta=\gamma-1 >\beta^*$.

We consider the figela process $(S_{0}(t),t\geq0)$ with autosimilarity
parameter \mbox{$\alpha=0$}.
We fix $t>0$ for the moment and denote
the maximal number of ends in a fragment of $S_0(t)$ by $\mathrm{E}(t)$.

Recall that $R_i^{(1,x)}(S_0(t))$, $1\leq i\leq H_{S_0(t)}(1,x)+1$
are the fragments of $S_0(t)$ separating $1$ from $x$.
Any chord
$[yz]$ with $(y,z)\in S^*(\infty)\cap(I\times J)$ must be contained
in the closure of one of these fragments
[otherwise this chord would cross one of the chords of $S_0(t)$, which
is impossible].
Consequently, the sets
\[
\bigl(I\cap\overline{R_i^{(1,x)}(S_0(t))}\bigr) \times\bigl(J\cap\overline
{R_i^{(1,x)}(S_0(t))}\bigr), \qquad
1\leq i\leq H_{S_0(t)}(1,x)+1
\]
form a covering of $S^*(\infty)\cap(I\times J)$. We get a finer
covering by considering the sets
$\overline C\times\overline D$, where $C$ varies over the connected
components of $I\cap{R_i^{(1,x)}(S_0(t))}$,
and $D$ varies over the connected components of $J\!\cap\!{R_i^{(1,x)}(S_0(t))}$.
We denote these connected components by $C_{ik}$, $1\leq k\leq k_i$ and
$D_{i\ell}$, \mbox{$1\leq\ell\leq\ell_i$}, respectively. Note that $k_i
\leq\mathrm{e}(R_i^{(1,x)}(S_0(t)))\leq\mathrm{E}(t)$, and the same
bound holds for $\ell_i$.
Summarizing the preceding discussion, we have
%
%
\begin{equation}
\label{covering1}
L_{(I,J)} \subset\Biggl( \bigcup_{i=1}^{H_{S_0(t)}(1,x)+1} \bigcup
_{k=1}^{k_i}\bigcup_{\ell=1}^{\ell_i}
{\mathcal C}^i_{k,\ell} \Biggr),
\end{equation}
where $ {\mathcal C}^i_{k,\ell}$ stands for the union of all chords
$[yz]$ for $y\in\overline C_{ik}$
and $z\in\overline D_{i\ell}$.

For every $1\leq i\leq H_{S_0(t)}(1,x)+1$, let
\[
\eta_i(t) =2\pi\mathrm{m}\bigl(R_i^{(1,x)}(S_0(t))\bigr)
\]
be the length of $R_i^{(1,x)}(S_0(t))\cap\sun$. Obviously the length
of any of the arcs~$C_{ik}$, $D_{i\ell}$ is bounded above by $\eta_i(t)$. Consequently, we
can cover each set
$ {\mathcal C}^i_{k,\ell}$ by at most $2 \eta_i(t)^{-1}$ disks of
diameter $2\eta_i(t)$.
From this observation and (\ref{covering1}), we get a covering of
$L_{(I,J)}$ by disks
of diameter at most $2\max\{\eta_i(t)\dvtx 1\leq i\leq
H_{S_0(t)}(1,x)+1\}
$, such that the sum
of the $\gamma$th powers of the diameters of disks in this covering
is bounded above by
%
%
\begin{equation}
\label{Haus-majo}
2^{1+\gamma}\mathrm{E}(t)^2 \sum_{i=1}^{H_{S_0(t)}(1,x)+1} \eta
_i(t)^{\beta}.
\end{equation}

We then need obtain a bound for $\mathrm{E}(t)$. In the genealogical tree
of fragments, the number of ends
of a given fragment is at most the number of ends of its ``parent''
plus $1$. Consequently $\mathrm{E}(t)$ is smaller than the largest
generation of a fragment of $S_{0}(t)$. In our case $\alpha=0$, the genealogy
of fragments is described by a standard Yule process (indeed, each
fragment gives birth to two new fragments
at rate $1$). Easy estimates show that $\mathrm{E}(t) \leq t^2$ for all
large enough $t$, almost surely.
On the other hand, Theorem \ref{tech}(iv) implies that
\[
\limsup_{t\to\infty} \Biggl(\exp(\kappa_{\nu_{D}}(\beta
)t )\sum_{i=1}^{H_{S_0(t)}(1,x)+1}
\eta_i(t)^\beta\Biggr)<\infty  \qquad\mbox{a.s.}
\]
Since $\beta>\beta^*$ we have $\kappa_{\nu_{D}}(\beta)>0$. From
the preceding display
and the bound $\mathrm{E}(t)\leq t^2$ for $t$ large, we now deduce that
the quantity (\ref{Haus-majo})
tends to $0$ as \mbox{$t\to\infty$}. The upper bound $\op{dim}
L_{(I,J)}\leq\gamma$ follows.
By (\ref{HausdorffTech}) we have also $\op{dim} L_\infty\leq\gamma
$ and since
$\gamma> \beta^*+1$ was arbitrary, we conclude that $\op{dim}
L_\infty\leq\beta^*+1$.

\subsection*{Lower bound} For $(I,J)\in{\cal I}$, let $A_{(I,J)}$ be the
set of
all $y\in I$ such that there exists $z\in J$ with $(y,z)\in S^*(\infty
)$. By \cite{LGP08}, Proposition~2.3(i), we have
\[
\op{dim}(L_\infty)\geq\op{dim}\bigl(A_{(I,J)}\bigr)+1\vadjust{\goodbreak}
\]
for every $(I,J)\in{\cal I}$ (\cite{LGP08} deals with hyperbolic
geodesics instead of chords, but the argument
is exactly the same). For any rational $\delta\in\,]0,1/4[$, set
$I_\delta=\{e^{2i\pi r}\dvtx\delta\leq r\leq\frac{1}{2}-\delta\}$.
Also set $J_0=\{e^{2i\pi r}\dvtx\frac{1}{2}\leq r\leq1\}$. We will prove
that almost surely,
for all $\delta$ sufficiently small, we have
%
%
\begin{equation}
\label{HausdorffLB}
\op{dim}\bigl(A_{(I_\delta,J_0)}\bigr)\geq\beta^*.
\end{equation}
The desired lower bound for $\op{dim}(L_\infty)$ will then
immediately follow.

In order to get the lower bound (\ref{HausdorffLB}), we construct a
suitable random measure on
$A^1_{(I_\delta,J_0)}$. We define a finite random measure $\mu_\delta
$ on $[\delta,\frac{1}{2}-\delta]$
by setting, for every $r,s\in[\delta,\frac{1}{2}-\delta]$ with
$r\leq s$,
\[
\mu_\delta([r,s]) = \min_{u\in[s,\fraca{1}{2}]} \mathscr
{M}_{\infty}(e^{2i\pi u})
- \min_{u\in[r,\fraca{1}{2}]} \mathscr{M}_{\infty}(e^{2i\pi u}).
\]
Clearly, if $r$ belongs to the topological support of $\mu_\delta$,
we have
\[
\mathscr{M}_{\infty}(e^{2i\pi r})= \min_{u\in[r,\fraca{1}{2}]}
\mathscr{M}_{\infty}(e^{2i\pi u}),
\]
and thus there exists $s\in[\frac{1}{2},1]$ such that
\[
\mathscr{M}_{\infty}(e^{2i\pi r}) = \mathscr{M}_{\infty}(e^{2i\pi s})
= \min_{u\in[r,s]} \mathscr{M}_{\infty}(e^{2i\pi u}).
\]
Therefore, with the notation of the previous section, we have
$r\overset{g_\infty}{\sim}s$
and also $r\overset{g_\infty}{\approx}s$ from the proof of Theorem
\ref{codingL}. It follows that $(e^{2i\pi r},e^{2i\pi s})\in
S^*(\infty)$ and $e^{2i\pi r}\in A_{(I_\delta,J_0)}$.

To summarize, if we denote the image of $\mu_\delta$ under the mapping
$r\longrightarrow e^{2i\pi r}$ by $\nu_\delta$, the measure $\nu
_\delta$ is supported on $A_{(I_\delta,J_0)}$.
From the H\"{o}lder continuity properties of the process $\mathscr
{M}_\infty$, we immediately
get that for every $\varepsilon>0$ there exists a (random) constant
$C_\varepsilon$
such that the $\nu_\delta$-measure of any ball is bounded above by
$C_\varepsilon$
times the $(\beta^*-\varepsilon)$th power of the diameter of the
ball. The
lower bound (\ref{HausdorffLB}) now follows from standard results
about Hausdorff
measures, provided that we know that $\nu_\delta$ is nonzero for
$\delta>0$ small, a.s.
However the total mass of $\nu_\delta$ clearly converges to $\mathscr
{M}_\infty(-1)>0$
as $\delta\to0$. This completes the proof. 

\begin{rek}
A simplified version of the preceding arguments gives the dimension of
the set of
all feet of nondegenerate chords of $S^*(\infty)$
\[
\op{dim}\bigl\{x\in\sun\dvtx \exists y\in\sun\setminus\{x\}\dvtx
(x,y)\in
S^*(\infty)\bigr\}= \beta^*  \qquad\mbox{a.s.}
\]
(Compare with Lemma \ref{nochordfrom1}.)
\end{rek}

\section{Convergence of discrete models}\label{sec7}

In this section, we prove Theorem~\ref{discreteapprox}. A key tool is
the maximality property in Theorem~\ref{max}. We will also
need the following geometric lemma, which considers laminations that
are ``nearly maximal.''\vadjust{\goodbreak}

\begin{lemma}\label{quasimaximal} Let $S$ be a figela and $\varepsilon
\in\,]0,1[$. Suppose that all fragments of $S$ have mass smaller than
$\varepsilon/2\pi$ and at most $3$ ends. Consider an
arbitrary lamination
\[
L=\bigcup_{i\in I} [x_iy_i],
\]
where the chords $[x_iy_i]$ do not cross. Suppose that the chords
of the figela $S$ belong to the collection $\{[x_iy_i]\dvtx i\in I\}$, and
in particular $L_S\subset L$. Then any chord $[x_iy_i]$, $i\in I$
lies within Hausdorff distance less than $\varepsilon$ from a chord of
the figela $S$.
\end{lemma}

\begin{figure}

\includegraphics{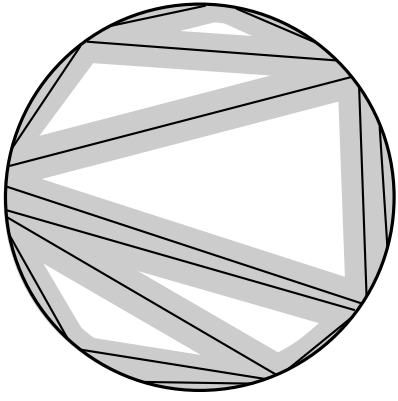}

\caption{Illlustration of the proof of lemma \protect\ref{quasimaximal}: The
chords $[x_iy_i]$ have to lie in the shaded
part of the figure.}
\label{fig6}
\end{figure}

We omit the easy proof, which should be clear from Figure~\ref{fig6}.

%
%

Let us turn to the proof of Theorem \ref{discreteapprox}. We fix
$\varepsilon>0$
and $\delta\in\,]0,1/2[$.

We use the genealogical structure of fragments as described in the
beginning of the proof of Theorem \ref{max}. We first observe that we
may fix an~inte\-ger~$m$ sufficiently large such that with probability at least
$1-\delta$ all the
fragments $R_u$ for $u\in\{0,1\}^m$ have mass less than $\varepsilon
/2\pi$. Then, using
Lemma~\ref{3bouts}, or rather the proof of this lemma, we can choose
an integer
$M\geq m$ large enough so that the following holds with probability
greater than
$1-\delta$: for every $u \in\{0,1\}^M$, there exists an integer
$j(u)\in\{m,\ldots,M\}$ such that the fragment $R_{[u]_{j(u)}}$ has
at most 3 ends.

From now on, we argue on the set where the preceding property holds and
where all the
fragments $R_u$ for $u\in\{0,1\}^m$ have mass less than $\varepsilon
/2\pi$.
For every $u\in\{0,1\}^M$, we choose the integer $j(u)$ as small as
possible and set $v(u)=[u]_{j(u)}$
to simplify notation. Then, if $u,u'\in\{0,1\}^M$, the fragments~$R_{v(u)}$ and $R_{v(u')}$
are either disjoint or equal. From this property, we easily get that
there exists a figela
$L^*$ whose fragments are the sets $R_{v(u)}$, $u\in\{0,1\}^M$. By
construction,
$L^*$ satisfies the assumptions of Lemma \ref{quasimaximal}.
Consequently, every chord
appearing in the figela process lies within distance at most
$\varepsilon$ from a chord of $L^*$.

For a given value of the integer $n\geq3$, consider now the discrete
triangulation $\Lambda_n$ of Section~\ref{sec1}, such that feet of
chords belong
to the set of $n$th roots of unity, and recall the recursive
construction of $\Lambda_n$. In this model we can introduce a
labelling of fragments analogous
to what we did in the continuous setting. For instance, $R^n_0$ and $R^n_1$
will be the fragments created by the first chord, ordered in a random
way. Then we look for the
first chord that falls in $R_0$ (if any) and call $R^n_{00}$ and
$R^n_{01}$ the new fragments
created by this chord, and so on. In this way we get a collection
$(R^n_u)_{u\in{\mathbb T}^{(n)}}$,
which is indexed by a random finite subtree ${\mathbb T}^{(n)}$ of
$\mathbb{T}$. It is easy to verify that, for every
integer $p\geq0$, ${\mathbb P}[\mathbb{T}_p\subset{\mathbb
T}^{(n)}]$ tends to $1$
as $n\to\infty$.

For every $u\in\mathbb{T}$, write $x_u$ and $y_u$ for the feet of the
first chord that will
split $R_u$ (again ordered in a random way). Introduce a similar
notation $x^n_u$ and $y^n_u$ in the discrete setting
[then of course $x^n_u$ and $y^n_u$ are only defined when $u\in\mathbb
{T}^{(n)}$
and $u$ is not a leaf of $\mathbb{T}^{(n)}$]. Since feet of chords are
chosen recursively
uniformly over possible choices, both in the discrete and in the continuous
setting, it should be clear that, for every integer $p\geq0$,
%
%
\begin{equation}
\label{convfeet}
((x^n_u,y^n_u) )_{u\in\mathbb{T}_p}
\overset{(d)}{\underset{n\to\infty}{\longrightarrow}}
((x_u,y_u) )_{u\in\mathbb{T}_p}.
\end{equation}
We apply this convergence with $p=M$. Using the Skorokhod
representation theorem,
we may assume that the preceding convergence holds almost surely. Then
almost surely for
$n$ sufficiently large, every chord of the figela $L^*$ (which must be
of the form
$[x_uy_u]$ for some $u\in\mathbb{T}_M$) lies within distance at most~$\varepsilon$
from a chord of $\Lambda_n$. Recalling the beginning of the proof, we
see that,
on an event of probability at least $1-2\delta$, every chord appearing
in the figela process lies within distance $2\varepsilon$ from a chord
of $\Lambda_n$, for all
$n$ sufficiently large.

We still need to prove the converse: We
argue on the same event of probability at least $1-2\delta$ and verify
that, for $n$ sufficiently large, every chord of~$\Lambda_n$ lies
within distance $2\varepsilon$ from the set $L_\infty$.
To this end, we use a symmetric argument. Assuming that $n$
is large enough so that $\mathbb{T}_M\subset{\mathbb T}^{(n)}$, we
let~$\Lambda^*_n$
be the figela whose fragments are the sets $R^n_{v(u)}$, $u\in\{0,1\}
^M$. The (almost sure)
convergence (\ref{convfeet}) guarantees that every chord of the figela $L^*$
is the limit as $n\to\infty$ of the corresponding chord of $\Lambda
^*_n$. It follows that,
for $n$ sufficiently large, $\Lambda^*_n$ satisfies the assumptions of
Lemma~\ref{quasimaximal},
and thus every chord of~$\Lambda_n$ lies within distance at most
$\varepsilon$
from a chord of $\Lambda^*_n$. Taking $n$ even larger if necessary, we
get that
every chord of~$\Lambda_n$ lies within
distance at most $2\varepsilon$ from a chord of $L^*$. This completes
the proof of the
first assertion of Theorem~\ref{discreteapprox}.

The second assertion is proved in a similar manner. Plainly, a
uniformly distributed random
permutation of $\{1,2,\ldots,n\}$ can be generated by first choosing
$\sigma(1)$ uniformly
over $\{1,\ldots,n\}$, then $\sigma(2)$ uniformly over $\{1,\ldots
,n\}\setminus\{\sigma(1)\}$, and so on.
From this simple remark, we see that the analogue of the convergence
(\ref{convfeet}) still holds for the feet of chords of the figela
$\widetilde\Lambda_n$. The remaining
part of the argument goes through without change.

\section{Extensions and comments}\vspace*{1pt}\label{sec8}

\subsection{\texorpdfstring{Case $\alpha=0$}{Case alpha = 0}} Recall from Theorem \ref{tech}(iv) the
definition of $\mathscr{H}_{0}(x)$ as the almost sure
limit of $e^{-t/3} H_{S_{0}(t)}(1,x) $ as $t\to\infty$. Note that
$\mathscr{H}_{0}(x)$ is an analogue
in the homogenous case $\alpha=0$ of $\mathscr{M}_{\infty}(x)$. In a
way similar to what we did for
$\mathscr{M}_{\infty}(x)$, one can verify that $\mathbb{E} [
\mathscr {H}_{0}(x)^p ]<\infty$ for every real $p\geq1$, and derive
integral equations for the moments
$h_{p}(r) = \mathbb{E} [ \mathscr{H}_{0}(e^{2i\pi r})^p ]$, for
$0\leq r\leq1$.
In the case $p=1$ we get
\[
\label{eqh1}
\frac{4}{3}h_{1}(r)=\int_{0}^r du \biggl(\frac
{1-r}{1-u} \biggr)^{2}\frac{2h_{1}(u)}{1-u}
+ \int_{r}^1 du \biggl(\frac{r}{u} \biggr)^{2}\frac{2h_{1}(u)}{u}.
\]
By differentiating this equation three times with respect to the
variable $r$ we get
\[
\frac{2}{3}h_{1}'''(r)= h_{1}''(r) \biggl(\frac{1}{1-r}-\frac
{1}{r} \biggr),
\]
leading to the explicit formula
\[
h_{1}(r) = \frac{8}{\pi}\sqrt{r(1-r)}.
\]
For higher values of $p$, we get the following bounds.
\begin{proposition}
\label{momentsbis} For every integer $p \geq1$ and every $\varepsilon
>0$, there exists a constant $K$ such that for every $r\in[0,1]$,
\[
h_{p}(r) \leq K \bigl(r(1-r)\bigr)^{\fracc{2p}{p+3}-\varepsilon}.
\]
\end{proposition}

We omit the proof, which uses arguments similar to the proof of
Proposition~\ref{moments}.
The bounds of Proposition \ref{momentsbis} are not sharp. Still they
are good enough to apply Kolmogorov's continuity criterion in order
to get a continuous modification of the process $(\mathscr
{H}_{0}(x))_{x\in\sun}$.

\subsection{Recursive self-similarity}

Set $Z_t=\mathscr{M}_\infty(e^{2i\pi t})$ for every $t\in[0,1]$. A~slightly more precise
version of Proposition \ref{chord} shows that the process $(Z_t)_{t\in
[0,1]}$ satisfies the following remarkable
self-similarity property. Let $Z'$ and $Z''$ be two independent copies
of $Z$ and let $(U_1,U_2)$ be distributed according
to the density $2{\mathbf1}_{\{0<u_1<u_2<1\}}$ and independent of
the pair $(Z,Z')$. Then the process $(\widetilde Z_t)_{t\in[0,1]}$
defined by
\[
\widetilde Z_t =
\cases{
\displaystyle
\bigl(1-(U_2-U_1)\bigr)^{\beta^*}Z'_{t/(1-(U_2-U_1))} ,\vspace*{1pt} \cr
\qquad \mbox{if }0 \leq t
\leq U_1, \vspace*{3pt}\cr\displaystyle
\bigl(1-(U_2-U_1)\bigr)^{\beta^*}Z'_{U_1/(1-(U_2-U_1))} + (U_2-U_1)^{\beta
^*} Z''_{(t-U_1)/(U_2-U_1)} ,\vspace*{1pt}\cr
\qquad \mbox{ if }U_1 \leq t \leq U_2,  \vspace*{3pt}\cr
\displaystyle
\bigl(1-(U_2-U_1)\bigr)^{\beta^*}Z'_{(t-(U_2-U_1))/(1-(U_2-U_1))} ,\vspace*{1pt}\cr
\qquad \mbox{ if }U_2 \leq t \leq1,
}
\]
has the same distribution as $(Z_t)_{t\in[0,1]}$.

Informally, this means that we can write a decomposition of $Z$ in two
pieces according to the following device. Throw two
independent uniform points $U_1$ and $U_2$ in $[0,1]$. Condition on the
event $U_1<U_2$ and set $M=1-(U_2-U_1)$. Then start from a (scaled)
copy of $Z$ of duration $[0,M]$ and ``insert'' at time $U_1$ another
independent scaled copy of $Z$ of duration $1-M$. Then the resulting
random function has the same distribution as $Z$.

In \cite{Ald94a}, Aldous describes such a decomposition in three
pieces for the Brownian excursion, which is closely related
to the random geodesic lamination~$L_\mathbf{e}$ of Theorem \ref
{AldousTh}. Aldous also
conjectures that there cannot exist a decomposition of the Brownian
excursion in two pieces of the type described above.

It would be interesting to know whether the preceding decomposition of~$Z$ (along with some regularity properties)
characterizes the distribution of~$Z$ up to trivial scaling constants.
One may also ask whether the scaling exponent
$\beta^*$ is the only one for which there can exist such a
decomposition in two pieces.

\section*{Acknowledgments} We are grateful to Gr\'{e}gory Miermont for a
number of
fruitful discussions. We also thank Jean Bertoin, Fran\c{c}ois David and
Kay J\"{o}rg Wiese for useful conversations. The second author is
indebted to
Fr\'{e}d\'{e}ric Paulin for suggesting the study of recursive
triangulations of the disk.
Finally we thank the referee for several useful comments.

%

\printaddresses

\end{document}